\documentclass[reqno, 11pt]{article}

\RequirePackage[OT1]{fontenc}
\RequirePackage{amsthm,amsmath}
\RequirePackage[numbers]{natbib}
\RequirePackage[colorlinks,citecolor=blue,urlcolor=blue]{hyperref}

\usepackage{mathrsfs}
\usepackage{amssymb}
\usepackage{amsfonts}
\usepackage{stmaryrd}
\usepackage{mathrsfs,dsfont}
\usepackage{bbm}

\usepackage{a4wide}
\def\ud{\mathrm{d}}

\def\naturals{\mathbb{N}}
\def\reals{\mathbb{R}}

\def\relatives{\mathbb{Z}}
\def\pp{\mathsf{P}}
\def\ee{\mathsf{E}}
\def\ff{\mathsf{F}}
\def\bb{\textsf{B}}
\def\gg{\textsf{G}}
\def\hh{\textsf{H}}
\def\dd{\textsf{D}}
\def\vv{\textsf{V}}
\def\probabilityspace{(\Omega, \mathscr{F}, \textsf{P})}

\def\pmspace{\mathcal{P}(0,1)}

\def\xsng{\overline{x}_{n,\gamma}}
\allowdisplaybreaks

\numberwithin{equation}{section}
\theoremstyle{plain}
\newtheorem{lem}{Lemma}[section]
\newtheorem{prp}{Proposition}[section]
\newtheorem{thm}{Theorem}[section]
\newtheorem{cor}{Corollary}[section]
\allowdisplaybreaks

\makeatletter
\newcommand{\subjclass}[2][2010]{%
  \let\@oldtitle\@title%
  \gdef\@title{\@oldtitle\footnotetext{#1 \emph{Mathematics subject classification.} #2}}%
}
\newcommand{\keywords}[1]{%
  \let\@@oldtitle\@title%
  \gdef\@title{\@@oldtitle\footnotetext{\emph{Key words and phrases.} #1.}}%
}
\makeatother

\author{Emanuele Dolera and Stefano Favaro}

\title{\textbf{Rates of convergence in de Finetti's representation theorem, and Hausdorff moment problem}}

\subjclass{60F05, 60G09}

 \keywords{de Finetti's law of large numbers, de Finetti's representation theorem, Edgeworth expansions, exchangeability, Hausdorff moment problem, Kolmogorov distance, Wasserstein distance.}

\date{}

\begin{document}
\maketitle

\begin{abstract}
Given a sequence $\{X_n\}_{n \geq 1}$ of exchangeable Bernoulli random variables, the celebrated de Finetti representation theorem states that $\frac{1}{n} \sum_{i = 1}^n X_i \stackrel{a.s.}{\longrightarrow} Y$ for a suitable random variable $Y : \Omega \rightarrow [0, 1]$ satisfying $\pp[X_1 = x_1, \dots, X_n = x_n\ |\ Y] = Y^{\sum_{i=1}^{n}x_{i}} (1 - Y)^{n - \sum_{i=1}^{n}x_{i}}$. In this paper we study the rate of convergence in law of $\frac{1}{n} \sum_{i = 1}^n X_i$ to $Y$ under the Kolmogorov distance. After showing that a rate of the type of $1/n^{\alpha}$ can be obtained for any index $\alpha \in (0,1]$, we find a sufficient condition on the distribution of $Y$ for the achievement of the optimal rate of convergence, that is $1/n$. Besides extending and strengthening recent results on the rate of convergence in de Finetti's representation theorem under the weaker Wasserstein distance, our main result weakens the regularity hypotheses on $Y$ in the context of the Hausdorff moment problem.
\end{abstract}

%


\section{Introduction}

This paper contributes to the study of the rate of convergence of the law of large numbers for exchangeable random variables (r.v.'s) in the sense of de Finetti \cite{deFin33}. For ease in exposition, we confine ourselves to an infinite sequence $\{X_n\}_{n \geq 1}$ of Bernoulli variables defined on the probability space $\probabilityspace$. The sequence 
$\{X_n\}_{n \geq 1}$ satisfies the exchangeability condition if there holds
$$
\pp[X_1 = x_1, \dots, X_n = x_n] = \pp[X_1 = x_{\sigma_n(1)}, \dots, X_n = x_{\sigma_n(n)}]
$$ 
for all $n \in \naturals$, $(x_1, \dots, x_n) \in \{0, 1\}^n$ and permutation $\sigma_n$ of the set $\{1, \dots, n\}$. De Finetti \cite{deFin33} proved a strong law of large numbers for the exchangeable $X_{i}$'s, i.e. $\frac{1}{n} \sum_{i = 1}^n X_i \stackrel{a.s.}{\longrightarrow} Y$ for a suitable  r.v. $Y : \Omega \rightarrow [0, 1]$ satisfying 
$\pp[X_1 = x_1, \dots, X_n = x_n\ |\ Y] = Y^{s_n} (1 - Y)^{n - s_n}$, where $s_n := \sum_{i = 1}^n x_i$. This identity yields the so-called de Finetti representation theorem \cite{deFin30}, which reads as
\begin{equation} \label{eq:deFin}
\pp[X_1 = x_1, \dots, X_n = x_n] = \int^1_0 \theta^{s_n} (1 - \theta)^{n - s_n} \mu(\ud\theta)
\end{equation}
for all $n \in \naturals$ and $(x_1, \dots, x_n) \in \{0, 1\}^n$, where $\mu$, the so-called de Finetti measure (or prior measure), stands for the probability distribution (p.d.) of $Y$. See Aldous \cite{ald} and references therein for a comprehensive treatment of exchangeability and de Finetti's theorem. The above law of large number entails that the p.d. of $\frac{1}{n} \sum_{i = 1}^n X_i$, say $\mu_n$, converges weakly to $\mu$ as $n\rightarrow+\infty$ ($\mu_n \Rightarrow \mu$ in symbols), meaning that $\lim_{n \rightarrow \infty} \int_0^1 \psi(\theta) \mu_n(\ud\theta) = \int_0^1 \psi(\theta) \mu(\ud\theta)$ is valid for all $\psi : [0,1] \rightarrow \reals$ bounded and continuous. 

The study of the rate of convergence of the empirical measure $\mu_n$ to $\mu$ requires the choice of a suitable distance that induces the weak convergence. A reasonable choice of such a distance yields an explicit evaluation of the discrepancy between $\mu_n$ and $\mu$, as a function of the sample size $n$, and also a practical interpretation of the approximation from various points of view. In this paper we focus on the Kolmogorov distance which, for any pair $(\nu_1, \nu_2)$ of probability measures (p.m.'s) on $([0,1], \mathscr{B}([0,1]))$, is defined as 
$$
\ud_K(\nu_1; \nu_2) :=  \sup_{x \in [0,1]} | \nu_1([0,x]) - \nu_2([0,x]) |  = \sup_{x \in [0,1]} | \ff_1(x) - \ff_2(x) |\ .
$$   
We recall that $\ud_K$ metrizes the weak convergence on the space $\pmspace$ of all p.m.'s on $([0,1], \mathscr{B}([0,1]))$ when the limiting p.m. has a continuous distribution function (d.f.). The next theorem is the main result of the present paper. It provides the first quantitative version, with respect to the Kolmogorov distance, of de Finetti's law of large numbers for $\mu_{n}$. For completeness, we denote by $\mathrm{L}^{\infty}(0,1)$ the space of essentially bounded function on $(0,1)$, with respect to the Lebesgue measure.

\begin{thm} \label{thm:peccati}
If $\mu$ has density function $f$ with distributional derivative $f'$ such that 
\begin{equation} \label{eq:peccatiH}
[\theta (1 - \theta)]^{\gamma} |f'(\theta)| \in \mathrm{L}^{\infty}(0,1),
\end{equation}
for some $\gamma \in (0,1)$, then
\begin{equation} \label{eq:MAIN} 
\ud_K(\mu_n; \mu) \leq \frac{C(\mu)}{n}
\end{equation}
is valid for all $n \in \naturals$, with $C(\mu)$ that depends on $\mu$ only through \emph{ess}.$\sup_{\theta \in [0,1]} f(\theta)$ and
\emph{ess}.$\sup_{\theta \in [0,1]} [\theta (1 - \theta)]^{\gamma} |f'(\theta)|$.
\end{thm}

Despite the long history of the celebrated de Finetti representation theorem, the study of the rate of convergence of $\mu_n$ to $\mu$ has been initiated very recently in the work of Mijoule, Peccati and Swan \cite{MPS}. They proved a quantitative version of de Finetti's law of large numbers with respect to the Kantorovich distance (or Wasserstein distance of order $1$) $\ud_W(\nu_1; \nu_2) := \int^1_0 | \ff_1(x) - \ff_2(x) | \ud x$. In particular, they showed that for any $n \in \naturals$
$$
\frac{C_1(\mu)}{n} \leq \ud_W(\mu_n; \mu) \leq \sqrt{\frac{C_1(\mu)}{n}}
$$ 
hold for any $\mu \in \pmspace$, where $C_1(\mu) := \int_0^1 \theta(1-\theta) \mu(\ud \theta)$. Thus, $1/n$ is the best possible rate of convergence to zero also for $\ud_K(\mu_n; \mu)$, since $\ud_K(\mu_n; \mu) \geq \ud_W(\mu_n; \mu)$, and it follows that $\ud_W(\mu_n; \mu)$ goes to zero at least as fast as $1/\sqrt{n}$. Mijoule Peccati and Swan \cite{MPS} also provided sufficient conditions for the achievement of the best rate of convergence. In particular, if the p.d. $\mu$ is absolutely continuous with a density $f$ satisfying $\int_0^1 \theta(1- \theta) |f'(\theta)| \ud\theta < +\infty$, where $f'$ stands for the (distributional) derivative of $f$, then 
$$
\ud_W(\mu_n; \mu) \leq \frac{C_2(\mu)}{n}
$$
holds for any $n \in \naturals$ with an explicit constant $C_2(\mu)$. Finally, it is proved that, for every $\delta \in [\frac{1}{2}, 1]$, there exists a suitable $\mu \in \pmspace$ for which $\ud_W(\mu_n; \mu) \sim 1/n^{\delta}$ 
as $n \rightarrow +\infty$. The work of Mijoule Peccati and Swan \cite{MPS} generalizes bound of the form $\ud_W(\mu_n; \mu) \leq C/n$, which was originally obtained in Goldstein and Reinert \cite{GoldRein} under the assumption that $\mu$ is a beta distribution. See also D\"obler \cite{Dobler} and references therein for related results. We recall that the beta distribution with parameters $(a, b) \in (0,+\infty)^2$ is the element of $\pmspace$ corresponding to the probability density function 
\begin{equation} \label{beta_density}
(0,1) \ni \theta \mapsto \beta(\theta; a, b) := \frac{\Gamma(a+b)}{\Gamma(a)\Gamma(b)} \theta^{a - 1} (1 - \theta)^{b - 1}\ .
\end{equation} 

Our study of $\ud_K(\mu_n; \mu)$ is more challenging than the study of $\ud_W(\mu_n; \mu)$. While results in Mijoule Peccati and Swan \cite{MPS} rely on classical Berry-Esseen bounds for the Gaussian approximation in the central limit theorem, the proof of Theorem \ref{thm:peccati} relies on novel, refined version of these Berry-Esseen bounds, usually known as Edgeworth expansions. See Chapter 5 and Chapter 6 in Petrov \cite{petrov75} or Chapter 3 in Ibragimov and Linnik \cite{IbraLinn}. To highlight the difference between the study of $\ud_K(\mu_n; \mu)$ and $\ud_W(\mu_n; \mu)$, we state a simple result on beta prior measures, whose proof can be obtain by direct computation. This shows that any rate $n^{-\alpha}$, with $\alpha \in (0,1]$, is actually achieved by $\ud_K(\mu_n; \mu)$. 

\begin{prp} \label{prop:beta}
If $\mu$ is the beta distribution with parameter $(\alpha, 1)$ or $(1, \alpha)$, with $\alpha > 0$, then there exists a constant $C_{\alpha}$ for which, for any $n \in \naturals$,
\begin{equation} \label{eq:beta}
\ud_K(\mu_n; \mu) \leq C_{\alpha} \left(\frac{1}{n}\right)^{\alpha \wedge 1} 
\end{equation}
is fulfilled, where $\wedge$ denotes the minimum value. 
\end{prp}

We stress that, in view of results in Goldstein and Reinert \cite{GoldRein} and Mijoule Peccati and Swan \cite{MPS}, for the beta prior we have $\ud_W(\mu_n; \mu) \sim 1/n$. This shows that the asymptotic behaviour of $\ud_K(\mu_n; \mu)$ is different from the asymptotic behaviour of $\ud_W(\mu_n; \mu)$. In addition, we notice that it seems not convenient at all to resort to the inequality 
\begin{equation} \label{eq:KminW} 
\ud_K(\nu_1; \nu_2) \leq C(\nu_2) \sqrt{\ud_W(\nu_1; \nu_2)}, 
\end{equation}
which is valid whenever $\nu_2$ has a bounded density. See, e.g., Gibbs and Su \cite{GibbsSu}. Indeed, for the beta distribution with parameters $(\alpha, 1)$ or $(1, \alpha)$, with $\alpha \geq 1$, \eqref{eq:KminW} would lead to an upper bound like $C/\sqrt{n}$, which is worse than the upper bound $C/n$ given by \eqref{eq:beta}.

Besides extending and strengthening results obtained in Goldstein and Reinert \cite{GoldRein} and Mijoule Peccati and Swan \cite{MPS}, Theorem \ref{thm:peccati} is also connected to the celebrated Hausdorff moment problem.  The Hausdorff moment problem is known to be closely related to de Finetti's theorem. See, e.g., Ressel \cite{Res(85)}, Diaconis and Freedman \cite{Dia(04a)},   Diaconis and Freedman \cite{Dia(04b)} and references therein. Within the context of the Hausdorff moment problem, the main result of  Mnatsakanov \cite{Mna} shows the rate of convergence of $\ud_K(\mu_n; \mu)$ under a certain assumption on the prior measure $\mu$. In the equivalent reformulation of our problem as the finding of the rate of approximation in the Hausdorff moment problem, Theorem 2 in Mnatsakanov \cite{Mna} provides the existence of another constant $C_*(\mu)$ 
for which \eqref{eq:MAIN} holds for any $n \in \naturals$, with another constant $C_*(\mu)$ in the place of our $C(\mu)$. In spite of a very direct proof, the main difference is that 
$C_*(\mu)$ depends on ess.$\sup_{\theta \in [0,1]} |f'(\theta)|$, which is tantamount to requiring that the density $f$ of $\mu$ belongs to $\mathrm{W}^{1,\infty}(0,1)$, the Sobolev space of essentially bounded functions on $(0,1)$ with an essentially bounded distributional derivative. The comparison of Theorem 2 in Mnatsakanov \cite{Mna} with Proposition \ref{prop:beta} shows that the assumption $f \in \mathrm{W}^{1,\infty}(0,1)$ is indeed too strong and far from capturing the whole class of prior distributions for which \eqref{eq:MAIN} is met. Theorem \ref{thm:peccati} fills this gap by providing the general sufficient condition \eqref{eq:peccatiH} for the achievement of the best rate $1/n$. This leads to a remarkable improvement of Theorem 2 in Mnatsakanov \cite{Mna}. For prior distributions $\mu$ with a support strictly contained in $(0,1)$, condition \eqref{eq:peccatiH} boils down to the assumption that  $f \in \mathrm{W}^{1,\infty}(0,1)$, but, without this restriction, it is evident that Theorem \ref{thm:peccati} improves Theorem 2 in Mnatsakanov \cite{Mna} by allowing $f'(\theta)$ to diverge moderately at $\theta = 0$ and $\theta = 1$.
As a final remark, we show that Mnatsakanov's result cannot be re-adapted via a smoothing argument to obtain our sharper bound, thus justifying the effort of providing a longer and more complex proof. In fact, given any prior $\mu$ with density $f \not\in \mathrm{W}^{1,\infty}(0,1)$, one could try to smooth it by introducing a family of new priors, say 
$\{\mu_{\varepsilon}\}_{\varepsilon > 0}$, each with density $f_{\varepsilon} \in \mathrm{W}^{1,\infty}(0,1)$, so that $\mu_{\varepsilon} \Rightarrow \mu$ as $\varepsilon \downarrow 0$. Letting $\{X_n^{(\varepsilon)}\}_{n \geq 1}$ be a new sequence of exchangeable Bernoulli r.v.'s having $\mu_{\varepsilon}$ as de Finetti's measure, one could argue by resorting to the following triangular inequality:
$$
\ud_K(\mu_n; \mu) \leq \ud_K(\mu_n; \mu_{n,\varepsilon}) + \ud_K(\mu_{n,\varepsilon}; \mu_{\varepsilon}) + \ud_K(\mu_{\varepsilon}; \mu) 
$$
where $\mu_{n,\varepsilon}(\cdot) := \pp\big[\frac 1n \sum_{i=1}^n X_i^{(\varepsilon)} \in \cdot \big]$. Since $\ud_K(\mu_{n,\varepsilon}; \mu_{\varepsilon}) \leq C_*(\mu_{\varepsilon})/n$ would follow from Theorem 2 in Mnatsakanov \cite{Mna}, the achievement of \eqref{eq:MAIN} through this line of reasoning
would entail ess.$\sup_{\theta \in [0,1]} |f_{\varepsilon}'(\theta)| \leq C < +\infty$ for all $\varepsilon > 0$. 
Moreover, to have $\ud_K(\mu_{\varepsilon}; \mu) \leq C/n$, the parameter $\varepsilon$ should depend on $n$, yielding  
$\varepsilon = \varepsilon(n) \downarrow 0$ as $n \rightarrow +\infty$. But now, the Ascoli-Arzel\`a theorem implies that $f \in \mathrm{W}^{1,\infty}(0,1)$, leading to a contradiction.

As a corollary, we state a result for $\mu$ being a beta distribution. This agrees with Proposition \ref{prop:beta}, capturing exactly the elements of this class of priors for which the bound \eqref{eq:MAIN} is valid. Indeed, since the beta density function belongs to $\mathrm{W}^{1,\infty}(0,1)$ if and only if $a,b \geq 2$, we improve the assumption by the following corollary.

\begin{cor} 
If $\mu$ is the beta distribution with parameters $(a, b)$, then \eqref{eq:MAIN} is fulfilled if and only if $a, b \geq 1$.
\end{cor}

To conclude we state a proposition that deals with a larger class of priors than the class considered in Theorem \ref{thm:peccati}. Specifically, we show that a rate of convergence for $\ud_K(\mu_n; \mu)$, although not sharp, can be obtained also for a non absolutely continuous prior $\mu$, provided that its d.f. $\ff$ is H\"older continuous. This happens, for instance, if $\ff$ coincides with the Cantor function. The determination of the sharp rate in this setting remains an interesting open problem.

\begin{prp} \label{prop:Holder}
If $\mu \in \pmspace$ has a $\gamma$-H\"older continuous d.f. for some $\gamma \in (0,1]$,  then there exists a suitable constant $L_{\gamma}(\mu)$ for which, for any $n \in \naturals$
\begin{equation} \label{eq:Holder}
\ud_K(\mu_n; \mu) \leq \frac{L_{\gamma}(\mu)}{n^{\gamma/2}}
\end{equation}
\end{prp}

The paper is structured as follows. In Section 2 we present the proof of Theorem \ref{thm:peccati}, which requires a few preliminary lemmas. Proofs of Proposition \ref{prop:beta} and Proposition \ref{prop:Holder} are deferred to the appendix.


\section{Proof of Theorem \ref{thm:peccati}}\label{sec2}

We start with some preliminary lemmas. First, a decomposition lemma for probability density functions which will be used to justify the introduction of the additional hypothesis $f(0) = f(1) = 0$ in the first part of the proof of Theorem \ref{thm:peccati}. As to notation, any relation involving the symbols $A_{\pm}$ and $f_{\pm}$ must be intended as a short-hand for the two analogous relations which hold with $A_+, f_+$ and with $A_-, f_-$, respectively, in place of $A_{\pm}, f_{\pm}$.

\begin{lem} \label{lm:decomposition}
Given a probability density function $f$ on $[0,1]$ which is expressed by a polynomial, then there exist three non-negative constants $A_{\infty}$, $A_{+}$, $A_{-}$ and three continuous probability density functions $f_{\infty}$, $f_{+}$, $f_{-}$ on $[0,1]$ such that:
\begin{enumerate}
\item[i)] $A_{\infty} \leq \|f\|_{\infty} := \sup_{\theta \in [0,1]} |f(\theta)|$ and $A_{\pm} \leq 1+\|f\|_{\infty}$;
\item[ii)] $f_{\infty} \in \mathrm{W}^{1,\infty}(0,1)$ with $A_{\infty}\|f_{\infty}\|_{\infty} \leq \|f\|_{\infty}$ and $A_{\infty}\|f'_{\infty}\|_{\infty} \leq 2\|f\|_{\infty}$;
\item[iii)] $f_{+}, f_{-} \in \mathrm{W}^{1,\infty}(0,1)$, $f_{+}(0) = f_{+}(1) = f_{-}(0) = f_{-}(1) = 0$, $A_{\pm}\|f_{\pm}\|_{\infty} \leq 2\|f\|_{\infty}$ and, for any $\gamma \in (0,1)$, 
\begin{equation} \label{eq:pm}
A_{\pm} \sup_{\theta \in [0,1]} [\theta(1-\theta)]^{\gamma}|f'_{\pm}(\theta)| \leq \frac{2}{4^{\gamma}} \|f\|_{\infty} + \sup_{\theta \in [0,1]} [\theta(1-\theta)]^{\gamma}|f'(\theta)|\ ; 
\end{equation}
\item[iv)] $f(\theta) = A_{\infty}f_{\infty}(\theta) + A_{+}f_{+}(\theta) - A_{-}f_{-}(\theta)$ for all $\theta \in [0,1]$.
\end{enumerate}
\end{lem}

\begin{proof}
If $f(0) = f(1) = 0$, the thesis is trivial. Otherwise, we put $A_{\infty} = \frac{f(1)+f(0)}{2}\leq \|f\|_{\infty}$ and $f_{\infty}(\theta) = \frac{[f(1)-f(0)]\theta + f(0)}{A_{\infty}}$, so that ii) holds trivially. Then, recalling that, for any $a,b \in \reals$, $a = b + (a-b)_{+} - (b-a)_{+}$, where $x_{+} := \max\{0, x\}$, we set $g_{+}(\theta) := (f(\theta) - A_{\infty}f_{\infty}(\theta))_{+}$ and $g_{-}(\theta) := (A_{\infty}f_{\infty}(\theta) - f(\theta))_{+}$. Thus, we put $A_{\pm} = \int_0^1 g_{\pm}(\theta)\ud\theta$ and $f_{\pm}(\theta) = g_{\pm}(\theta)/A_{\pm}$ with the proviso that, if $A_{+} = 0$ ($A_{-} = 0$, respectively), the definition of $f_{+}$ ($f_{-}$, respectively) is arbitrary and can be chosen equal to $6x(1-x)$. By definition, point iv) is met along with $f_{+}(0) = f_{+}(1) = f_{-}(0) = f_{-}(1) = 0$. Moreover, we have $A_{\pm} \leq 1+A_{\infty}$ and point i) follows. To prove that $f_{+}, f_{-} \in \mathrm{W}^{1,\infty}(0,1)$, it is enough to notice that $f(\theta) = A_{\infty}f_{\infty}(\theta)$ for finitely many $\theta$'s, by virtue of the fundamental theorem of algebra, and, in the complement of this set, both $f_{+}$ and $f_{-}$ are again two polynomials. To check that $A_{\pm}\|f_{\pm}\|_{\infty} \leq 2\|f\|_{\infty}$, it is enough to observe that $\|g_{\pm}\|_{\infty} \leq \|f\|_{\infty} + A_{\infty}\|f_{\infty}\|_{\infty}$. Finally, we just note that, except on the finite set $\{\theta \in [0,1]\ |\ f(\theta) = A_{\infty}f_{\infty}(\theta)\}$, we have $A_{\pm} |f'_{\pm}(\theta)| \leq |f'(\theta)| + A_{\infty}|f'_{\infty}(\theta)|$, so that we can deduce the validity of \eqref{eq:pm} from the previous bounds.
\end{proof}


Another preliminary result deals with further regularity properties of the densities that satisfy \eqref{eq:peccatiH}. We observe that, since $[x(1-x)]^{-\gamma} \in \mathrm{L}^1(0,1)$ if $\gamma \in (0,1)$, the validity of \eqref{eq:peccatiH} entails 
$f \in \mathrm{W}^{1,1}(0,1)$ and, hence, the existence of a continuous version of the same density on the whole set $[0,1]$. See, e.g., Theorem 8.2 in Brezis \cite{brezis}.

\begin{lem} \label{lm:regularity}
Let $f$ be a probability density function satisfying \eqref{eq:peccatiH} for some $\gamma \in (0,1)$, and let $\ff(x) = \int_0^x f(y)\ud y$.  
Then, there exists a positive constant $R(\gamma)$ such that
\begin{equation} \label{eq:der2}
\sup_{\substack{x \in (0,1), \\ 0 < w < x(1-x)}} x(1-x) \Big{|} \frac{\ff(x+w) -2\ff(x) + \ff(x-w)}{w^2}\Big{|} \leq R(\gamma) |f|_{1,\gamma}
\end{equation}
is fulfilled with $|f|_{1,\gamma} :=$ \emph{ess}.$\sup_{\theta \in [0,1]} [\theta(1-\theta)]^{\gamma} |f'(\theta)|$. Moreover, if the additional condition $f(0) = f(1) = 0$ (referred to the continuous representative of $f$) is in force, then:
\begin{enumerate}
\item[i)] $f(\theta) \leq M(f) \theta^{1-\gamma}$,  $f(\theta) \leq M(f) (1-\theta)^{1-\gamma}$ hold for all $\theta \in [0,1]$, with $M(f) := \frac{2^{\gamma}}{1-\gamma}|f|_{1,\gamma} + 2^{1-\gamma}\|f\|_{\infty}$;
\item[ii)] $\ff(x) \leq M(f) \dfrac{x^{2-\gamma}}{2-\gamma}$,  $1 - \ff(x) \leq M(f) \dfrac{(1-x)^{2-\gamma}}{2-\gamma}$ hold for all $x \in [0,1]$.
\end{enumerate}
\end{lem}
\begin{proof}
Since $f \in \mathrm{W}^{1,1}(0,1)$ by virtue of \eqref{eq:peccatiH}, the Taylor formula with integral remainder can be applied to obtain
$$
\ff(x+w) -2\ff(x) + \ff(x-w) = \int_x^{x+w} (x+w-t) f'(t)\ud t + \int_x^{x-w} (x-w-t) f'(t)\ud t
$$
for all $x \in (0,1)$ and $w$ satisfying $0 < w < x(1-x)$. Whence,
\begin{equation} \label{eq:TaylorTaylor}
|\ff(x+w) -2\ff(x) + \ff(x-w)| \leq \int_x^{x+w} (w+t-x) |f'(t)|\ud t + \int_{x-w}^x (w+x-t) |f'(t)|\ud t\ .
\end{equation}
At this stage, we show explicitly how to bound the former integral when $x \in (0,1/2]$, the other cases being analogous. Since $0 < x < x+w < 3/4$, then we get $|f'(t)| \leq 4^{\gamma} |f|_{1,\gamma} t^{-\gamma}$ for all 
$t \in [x, x+w]$, leading to
\begin{gather} 
\int_x^{x+w} (w+t-x) |f'(t)|\ud t \nonumber \\
\leq 4^{\gamma} |f|_{1,\gamma} \left[ w \frac{(x+w)^{1-\gamma}-x^{1-\gamma}}{1-\gamma} + \frac{(x+w)^{2-\gamma}-x^{2-\gamma}}{2-\gamma} - x\frac{(x+w)^{1-\gamma}-x^{1-\gamma}}{1-\gamma}\right]\ . \nonumber
\end{gather}
Then, we put $\eta := w/x$ and we observe that $\eta \in (0,1/2)$, so that the expression inside the brackets can be written as
\begin{align} \label{eq:wxgamma}
&\notag w x^{1-\gamma} \frac{(1+\eta)^{1-\gamma}-1}{1-\gamma}\\
&\quad+ x^{2-\gamma}\left[ \frac{(1+\eta)^{2-\gamma}-1-(2-\gamma)\eta}{2-\gamma} - \frac{(1+\eta)^{1-\gamma}-1-(1-\gamma)\eta}{1-\gamma} \right]\ .
\end{align}
We conclude this argument by noticing that, for any $\eta$ satisfying $|\eta| \leq 1/2$ and any $\alpha > 0$ there exists a constant $H(\alpha)$ such that $|(1+\eta)^{\alpha}-1-\alpha\eta| \leq H(\alpha) \eta^2$. This remark implies that the expression in \eqref{eq:wxgamma} is bounded by the following quantity
\begin{displaymath}
\left[1 + \frac{H(2-\gamma)}{2-\gamma} + \frac{3H(1-\gamma)}{2(1-\gamma)}\right]w^2x^{-\gamma},
\end{displaymath}
yielding
\begin{align*}
&\sup_{\substack{x \in (0,1/2], \\ 0 < w < x(1-x)}} \frac{x(1-x)}{w^2} \int_x^{x+w} (w+t-x) |f'(t)|\ud t\\
&\quad\leq 2^{3\gamma-2}\left[1 + \frac{H(2-\gamma)}{2-\gamma} + \frac{3H(1-\gamma)}{2(1-\gamma)}\right] |f|_{1,\gamma}\ .
\end{align*}
As recalled, the treatment of the latter integral on the right-hand side of \eqref{eq:TaylorTaylor} for $x \in (0,1/2]$ is analogous. Lastly, when $x \in [1/2,1)$, it is enough to change the variable $t = 1-s$, obtaining
$\int_x^{x+w} (w+t-x) |f'(t)|\ud t = \int_{1-x-w}^{1-x} (w+1-x-s) |f'(1-s)|\ud s$ and $\int_{x-w}^x (w+x-t) |f'(t)|\ud t =\int_{1-x}^{1-x+w} (w+s -1+x) |f'(1-s)|\ud s$, where the integrals in the new variable $s$ are exactly the integrals studied above. To prove i), we just write $f(\theta) = \int_0^{\theta} f'(y)\ud y$ and, confining to the case that $\theta \in [0,1/2]$, we exploit \eqref{eq:peccatiH} in the form $|f'(\theta)| \leq 2^{\gamma} |f|_{1,\gamma} \theta^{-\gamma}$. This proves the first bound, after noticing that $f(\theta) \leq 2^{1-\gamma}\|f\|_{\infty}\theta^{1-\gamma}$ is valid for any $\theta \in [1/2,1]$. For the latter bound, we start from $-f(\theta) = \int_{\theta}^1 f'(y)\ud y$ and we argue in an analogous way. To prove ii), it is enough to integrate the bounds obtained in point i).
\end{proof}


The last preliminary result provides with a refinement of the well-known estimates of Berry-Esseen type for the characteristic function of a normalized sum of i.i.d., centered r.v.'s. In fact, the following statement can be viewed as a generalization of 
Lemma 4 in Chapter VI of Petrov \cite{petrov75} and Theorem 3.2.1(2) in Chapter 3 of Ibragimov and Linnik \cite{IbraLinn} in the case that the summands possess the $3+\delta$ absolute moment for some $\delta \in (0,1)$ but, in general, not the fourth moment. 
\begin{lem} \label{lm:BE}
Let $\{V_n\}_{n \geq 1}$ be a sequence of i.i.d. r.v's defined on $(\Omega, \mathscr{F}, \textsf{P})$ such that $\beta_{3+\delta} := \ee[|V_1|^{3+\delta}] < +\infty$ holds for some $\delta \in (0,1)$, along with $\ee[V_1] = 0$ and 
$\ee[V_1^2] =: \sigma^2 > 0$. Upon putting $\alpha_3 := \ee[V_1^3]$ and $\psi_n(\xi) := \ee\left[\exp\left\{i\xi\left(\sum_{k=1}^n V_k\right)/\sqrt{n\sigma^2}\right\}\right]$, there holds
\begin{equation} \label{eq:BE3delta}
\Big|\psi_n(\xi) - e^{-\xi^2/2}\Big\{1 + \frac{\alpha_3}{6\sqrt{n}\sigma^3}(i\xi)^3 \Big\}\Big| \leq Q(\delta)\frac{\beta_{3+\delta}}{n^{(1+\delta)/2}\sigma^{3+\delta}} |\xi|^{3+\delta}(1 + |\xi|^4) e^{-\xi^2/4}
\end{equation}
for any $\xi$ satisfying $|\xi| \leq \frac{1}{4} \sqrt{n}\left(\frac{\sigma^{3+\delta}}{\beta_{3+\delta}}\right)^{1/(1+\delta)}$\!\!\!\!, where $Q(\delta)$ is a numerical constant independent of $\xi$ and the p.d. of $V_1$.
\end{lem}

\begin{proof}
The proof is based on the arguments used to prove Lemma A.2 in Dolera, Gabetta and Regazzini \cite{dgr} and Lemma 3.1 in Dolera and Regazzini \cite{dore}. First, we put $\psi(\xi) := \ee\left[e^{i\xi V_1}\right]$ and we observe that $\psi(\xi) = 1 -\frac{\sigma^2}{2}\xi^2 + \frac{\alpha_3}{6}(i\xi)^3 + \rho_{\delta}(\xi)$, where
$$
|\rho_{\delta}(\xi)| \leq \frac{2^{1-\delta}\beta_{3+\delta}}{(1+\delta)(2+\delta)(3+\delta)}|\xi|^{3+\delta}\ .
$$
See, e.g., Theorem 1 in Section 8.4 of Chow and Teicher \cite{chte}. Whence,
$$
\psi\left(\frac{\xi}{\sqrt{n\sigma^2}}\right) = 1 -\frac{1}{2n}\xi^2 + \frac{\alpha_3}{6\sigma^3n^{3/2}}(i\xi)^3 + \rho_{n,\delta}(\xi)
$$
with
$$
|\rho_{n,\delta}(\xi)| \leq \frac{2^{1-\delta}\beta_{3+\delta}}{(1+\delta)(2+\delta)(3+\delta)\sigma^{3+\delta}n^{(3+\delta)/2}}|\xi|^{3+\delta}\ .
$$
Now, we notice that Lyapunov's inequality entails $\sigma^{3+\delta} \leq \beta_{3+\delta}$, while H\"older's inequality shows that
$$
\frac{|\alpha_3|}{\sigma^3}\left(\frac{\sigma^{3+\delta}}{\beta_{3+\delta}}\right)^{3/(1+\delta)} \leq\ \frac{\beta_3}{\sigma^3}\left(\frac{\sigma^{3+\delta}}{\beta_{3+\delta}}\right)^{3/(1+\delta)} \leq\
\frac{\beta_3}{\sigma^3}\left(\frac{\sigma^{3+\delta}}{\beta_{3+\delta}}\right)^{1/(1+\delta)} \leq\ 1
$$
where $\beta_3 := \ee[|V_1|^3]$. Therefore, for any $\xi$ satisfying $|\xi| \leq \frac{1}{4}\sqrt{n}\left(\frac{\sigma^{3+\delta}}{\beta_{3+\delta}}\right)^{1/(1+\delta)}$\!\!\!\!, we have
$$
\Big|-\frac{1}{2n}\xi^2 + \frac{\alpha_3}{6\sigma^3n^{3/2}}(i\xi)^3 + \rho_{n,\delta}(\xi) \Big| \leq \frac{5}{128} \ .
$$
Thanks to this bound, we are allowed to consider the principal logarithm $\mathrm{Log}(1+z) := -\sum_{k=1}^{\infty} \frac{(-z)^k}{k}$, $|z| < 1$ and,
since $\psi_n(\xi) = \left[\psi\left(\frac{\xi}{\sqrt{n\sigma^2}}\right)\right]^n$, we have:
\begin{align*}
\psi_n(\xi)&= \exp\left\{n\mathrm{Log}\left[\psi\left(\frac{\xi}{\sqrt{n\sigma^2}}\right)\right]\right\}\\
& = \exp\left\{n\mathrm{Log}\left[1 -\frac{1}{2n}\xi^2 + \frac{\alpha_3}{6\sigma^3n^{3/2}}(i\xi)^3 + \rho_{n,\delta}(\xi)\right]\right\} \nonumber \\
&=e^{-\xi^2/2} \exp\left\{\frac{\alpha_3}{6\sqrt{n}\sigma^3}(i\xi)^3\right\} e^{\tau_{n,\delta}(\xi)} \nonumber
\end{align*}
where $\tau_{n,\delta}(\xi)$ is defined to be
$$
n\rho_{n,\delta}(\xi) + n\left[-\frac{1}{2n}\xi^2 + \frac{\alpha_3}{6\sigma^3n^{3/2}}(i\xi)^3 + \rho_{n,\delta}(\xi)\right]^2 \Upsilon\left(-\frac{1}{2n}\xi^2 + \frac{\alpha_3}{6\sigma^3n^{3/2}}(i\xi)^3 + \rho_{n,\delta}(\xi)\right)
$$
with $\Upsilon(z) := -\sum_{k=2}^{\infty} \frac{(-z)^{k-2}}{k}$ for $|z| < 1$. At this stage, we put $u_3(\xi) := \frac{\alpha_3}{6\sqrt{n}\sigma^3}(i\xi)^3$ and $\Theta(z) := e^z-1-z$, and we exploit the elementary inequality $|e^z - 1| \leq |z|e^{|z|}$ 
to obtain
\begin{align*}
&\Big|\psi_n(\xi) - e^{-\xi^2/2}\Big\{1 + \frac{\alpha_3}{6\sqrt{n}\sigma^3}(i\xi)^3 \Big\}\Big|\\
&\quad\leq e^{-\xi^2/2} e^{|u_3(\xi)|} \big|e^{\tau_{n,\delta}(\xi)} - 1\big| + e^{-\xi^2/2}\big|\Theta(u_3)\big| \nonumber \\
&\quad\leq e^{-\xi^2/2} e^{|u_3(\xi)|+|\tau_{n,\delta}(\xi)|} |\tau_{n,\delta}(\xi)| + e^{-\xi^2/2}\Theta(|u_3|) \ . \nonumber
\end{align*}
To conclude, it is enough to notice that, for any $\xi$ satisfying $|\xi| \leq \frac{1}{4}\sqrt{n} \left(\frac{\sigma^{3+\delta}}{\beta_{3+\delta}}\right)^{1/(1+\delta)}$\!\!\!\!, we have $|u_3(\xi)| \leq \frac{1}{384}$,
$|\tau_{n,\delta}(\xi)| \leq \frac{1}{4}\xi^2$,
$$
|\tau_{n,\delta}(\xi)| \leq Q_1(\delta) \frac{\beta_{3+\delta}}{n^{(1+\delta)/2}\sigma^{3+\delta}} |\xi|^{3+\delta}(1 + |\xi|^4)
$$
and
$$
\Theta(|u_3|) \leq Q_2(\delta) \frac{\beta_{3+\delta}}{n^{(1+\delta)/2}\sigma^{3+\delta}} |\xi|^{3(1+\delta)}
$$
for suitable constants $Q_1(\delta)$ and $Q_2(\delta)$ independent of $\xi$ and the p.d. of $V_1$.
\end{proof}


The way is now paved for the study of Theorem \ref{thm:peccati}. The first part of the proof, which requires the major effort, is devoted to proving \eqref{eq:MAIN} when the density $f$ of $\mu$, in addition to \eqref{eq:peccatiH}, satisfies $f(0) = f(1) = 0$. We recall again that \eqref{eq:peccatiH} entails 
$f \in \mathrm{W}^{1,1}(0,1)$ and, hence, the existence of a continuous version of this density on the whole set $[0,1]$, by virtue of Theorem 8.2 in Brezis \cite{brezis}. Obviously, the additional assumption $f(0) = f(1) = 0$ is referred to this version. After these preliminaries, we get into the real proof by defining $I(n,\gamma) := [\xsng, 1 - \xsng]$, with $\xsng := (1/n)^{\frac{1}{2-\gamma}}$, which is a proper interval provided that $n \geq 4$. Then, after denoting by $\ff_n$ ($\ff$, respectively) the d.f. associated to $\mu_n$ ($\mu$, respectively), we split the original quantity as follows:
\begin{align}
\ud_K(\mu_n; \mu) &\leq \sup_{x \in [0,\xsng]} | \ff_n(x) - \ff(x) | + \sup_{x \in I(n,\gamma)} | \ff_n(x) - \ff(x) |\nonumber \\
&\quad\quad + \sup_{x \in [1 - \xsng,1]} | \ff_n(x) - \ff(x) | \nonumber \\
&\leq \ff_n(\xsng) + \ff(\xsng) + \sup_{x \in I(n,\gamma)} | \ff_n(x) - \ff(x) | \nonumber \\
&\quad+ [1 - \ff_n(1-\xsng)] + [1 - \ff(1-\xsng)]\ . \ \ \ \ \ \label{eq:splitfirst}
\end{align}
To bound $\ff(\xsng)$ and $[1 - \ff(1-\xsng)]$, we use point ii) of Lemma \ref{lm:regularity}, which gives: 
\begin{equation} \label{eq:extremeff}
\ff(\xsng) + [1 - \ff(1-\xsng)] \leq \frac{2M(f)}{2-\gamma} \cdot \frac{1}{n}
\end{equation}
for all $n \geq 4$. To bound $\ff_n(\xsng)$ and $1 - \ff_n(1-\xsng)$, we invoke equation \eqref{eq:feller} in the proof of Proposition \ref{prop:beta} which, in combination with point ii) of Lemma \ref{lm:regularity}, yields
\begin{align*}
\ff_n(\xsng) &\leq M(f) \int^1_0 \beta(y; \lfloor n\xsng \rfloor + 1, n - \lfloor n\xsng \rfloor) y^{2-\gamma} \ud y \\
&= M(f) \frac{\Gamma(n+1)}{\Gamma( \lfloor n\xsng \rfloor + 1) \Gamma(n - \lfloor n\xsng\rfloor)} \frac{\Gamma( \lfloor n\xsng \rfloor + 3-\gamma) \Gamma(n - \lfloor n\xsng \rfloor)}{\Gamma(n+3-\gamma)}
\end{align*}
where $\beta$ is the same as in \eqref{beta_density} and $\lfloor \cdot \rfloor$ denotes the integral part.
The last expression can be majorized by means of \emph{Wendel's inequalities} (see (5) in Qi and Luo \cite{Qi}) as
$$
M(f) \frac{\lfloor n\xsng\rfloor + 2 - \gamma}{n+2-\gamma} \left( \frac{\lfloor n\xsng \rfloor + 1}{n + 1} \right)^{1-\gamma} \left( \frac{n+2-\gamma}{n + 1} \right)^{\gamma}
$$
which, for all $n \geq 4$, is less than $\frac{12}{5} M(f)/n$. To study $1 - \ff_n(1-\xsng)$, we argue as in the proof of Lemma \ref{lm:regularity} by considering
the exchangeable sequence $\{\overline{X}_n\}_{n \geq 1}$, where $\overline{X}_n := 1 - X_n$ for all $n \in \naturals$. Since the de Finetti measure of this new sequence is the element of $\pmspace$ associated to the d.f. 
$1 - \ff(1-x)$, we resort again to \eqref{eq:feller} to obtain
\begin{align*}
1 - \ff_n(1-\xsng)& \leq \pp\Big[\frac{1}{n}\sum_{i=1}^n \overline{X}_i \leq \xsng \Big]\\
& = \int^1_0 \beta(y; \lfloor n\xsng \rfloor + 1, n - \lfloor n\xsng \rfloor) [1 - \ff(1-y)] \ud y\ .
\end{align*}
Thus, by using the latter bound stated in point ii) of Lemma \ref{lm:regularity} and arguing exactly as above, we conclude that $1 - \ff_n(1-\xsng) \leq \frac{12}{5} M(f)/n$. 


Now, we study $\sup_{x \in I(n,\gamma)} | \ff_n(x) - \ff(x) |$. First, we get $\pp[\sum_{i=1}^n X_i = k\ |\ Y = \theta] = {n\choose k} \theta^k (1-\theta)^{n-k}$ for any $k \in \{0, \dots, n\}$, thanks to de Finetti's representation. Since
$$
\pp\left[\sum_{i=1}^n X_i \leq nx\ \Big|\ Y = \theta\right] = \pp\left[\frac{\sum_{i=1}^n (X_i - \theta)}{\sqrt{n\theta(1-\theta)}} \leq \frac{n(x - \theta)}{\sqrt{n\theta(1-\theta)}}\ \Big{|}\ Y = \theta\right]\ ,
$$
we put 
\begin{displaymath}
u := u(x, \theta, n) := \frac{n(x - \theta)}{\sqrt{n\theta(1-\theta)}}
\end{displaymath}
and
\begin{displaymath}
\bb_n(y ;\theta) := \pp\left[\sum_{i=1}^n (X_i - \theta) \leq y\sqrt{n\theta(1-\theta)}\ \Big{|}\ Y = \theta\right]
\end{displaymath}
for $y \in \reals$. To study $\bb_n(\cdot ;\theta)$, we make the key remark that it coincides with the d.f. of a normalized sum of i.i.d., centered r.v.'s, so that we can employ well-known results pertinent to the central limit theorem, as stated in 
Chapter 8 of Gnedenko and Kolmogorov,\cite{GneKol}, Chapters 5-6 of Petrov \cite{petrov75}, Chapter 3 of Ibragimov and Linnik \cite{IbraLinn}, and in Osipov \cite{osipov}. In particular, mimicking the main theorem in Osipov \cite{osipov}, we introduce the functions $\gg_n(y; \theta) := \Phi(y) + \hh_n(y; \theta)$ and
\begin{eqnarray} 
\hh_n(y; \theta) &:=& \frac{1}{\sqrt{2\pi n\theta(1-\theta)}} e^{-\frac{1}{2}y^2} \Big{\{} \frac{1}{6}(1-2\theta)(1-y^2) \nonumber \\ 
&+& S(n\theta + y \sqrt{n \theta(1 - \theta)} ) \Big{[} 1 + \frac{1-2\theta}{6\sqrt{n\theta(1 - \theta)}}(y^3 - 3y) \Big{]}\Big{\}}\ , \label{eq:IL339}
\end{eqnarray}
where $\Phi(y) := \int_{-\infty}^y \frac{1}{\sqrt{2\pi}} e^{-x^2/2} \ud x$ and $S(x) := \lfloor x \rfloor - x + \frac{1}{2}$. Now, Theorem 2 in Chapter 8 of Gnedenko and Kolmogorov \cite{GneKol} (see also Theorem 2b in Chapter II of Esseen \cite{esseen}) 
entails $\sup_{y \in \reals} |\bb_n(y ;\theta) - \gg_n(y; \theta)| \leq \epsilon_n(\theta)$ 
and provides the existence of three numerical constants $\lambda_1, \lambda_2 > 0$ and $n_0 \in \naturals$ (independent of $n$ and $\theta$) such that 
\begin{eqnarray} 
\epsilon_n(\theta) &=& \lambda_1 \int_{-n}^n \Big{|} \frac{\hat{\bb}_n(\xi; \theta) - \hat{\gg}_n(\xi; \theta)}{\xi} \Big{|} \ud\xi + \frac{\lambda_2}{n} \sup_{y \not\in \mathcal{Y}(n,\theta)} \big{|} \frac{\partial}{\partial y} \gg_n(y; \theta)\big{|} \nonumber \\
&\leq& \lambda_1 \int_{-n}^n \Big{|} \frac{\hat{\bb}_n(\xi; \theta) - \hat{\gg}_n(\xi; \theta)}{\xi} \Big{|} \ud\xi + \frac{\lambda_2 G}{n \theta(1-\theta)} \ \ \ \ \ \ (n \geq n_0) \label{eq:epsn}
\end{eqnarray}
where: $\hat{\bb}_n(\xi; \theta)$ and $\hat{\gg}_n(\xi; \theta)$ are Fourier-Stieltjes transforms of $\bb_n(\cdot; \theta)$ and $\gg_n(\cdot; \theta)$, i.e., $\int_{-\infty}^{\infty} e^{i\xi y} \ud_y \bb_n(y; \theta)$ and $\int_{-\infty}^{\infty} e^{i\xi y} \ud_y \gg_n(y; \theta)$, respectively; $\mathcal{Y}(n,\theta) := \big\{\frac{k - n\theta}{\sqrt{n\theta(1-\theta)}}\ |\ k \in \relatives\big\}$ is the set of the discontinuities of both $\bb_n(y; \theta)$ and $\gg_n(y; \theta)$; $G>0$ is another constant (independent of $n$ and $\theta$). 
Thus, since $\ff_n(x) = \int_0^1 \bb_n(u(x, \theta, n) ;\theta) f(\theta)\ud\theta$, we write 
\begin{align}
\sup_{x \in I(n,\gamma)} |\ff_n(x) - \ff(x)| &\leq \sup_{x \in I(n,\gamma)} \Big{|} \int_0^1 \Phi(u(x, \theta, n)) f(\theta)\ud\theta - \ff(x)\Big{|}  \nonumber \\
&\quad+ \sup_{x \in I(n,\gamma)} \int_0^1 |\hh_n(u(x, \theta, n); \theta)| f(\theta)\ud\theta + \int_0^1 \epsilon_n(\theta) f(\theta)\ud\theta\ \ \ \ \ \label{eq:splitmain}
\end{align}
and we try to bound each term on the right-hand side. Apropos of the first term on the right-hand side of \eqref{eq:splitmain}, we introduce a Gaussian r.v. $Z_n : \Omega \rightarrow \reals$ with zero mean and variance $1/n$, independent of $Y$, so that
\begin{align*}
\int_0^1 \Phi(u(x, \theta, n)) f(\theta)\ud\theta &= \ee\left[\pp[Y + Z_n\sqrt{Y(1-Y)} \leq x\ |\ Y] \right]\\
& = \pp[Y + Z_n\sqrt{Y(1-Y)} \leq x]\ .
\end{align*}
This d.f. (in the $x$ variable) plays an important role also in Mijoule Peccati and Swan \cite{MPS}, where its closeness to $\ff$ is proved with respect to the Kantorovich distance (see Proposition 4.1 therein). In any case, the proof in Mijoule Peccati and Swan \cite{MPS} is strongly based on a dual representation of $\ud_W$, which does not have any analog for $\ud_K$. Therefore, we tackle the problem by a direct computation which, after exchanging the order of conditioning in the above identity and using some elementary algebra, leads to
$$
\pp[Y + Z_n\sqrt{Y(1-Y)} \leq x] = \int_0^{+\infty} \sqrt{\frac{n}{2\pi}} \exp\Big\{-\frac{n}{2}z^2\Big\} \cdot [\ff(\theta_1(x,z)) + \ff(\theta_2(x,z))] \ud z
$$
where 
$$
\theta_1(x,z) := \frac{2x + z^2 -z\sqrt{z^2 + 4x(1-x)}}{2(z^2+1)}
$$
and
$$
\theta_2(x,z) := \frac{2x + z^2 +z\sqrt{z^2 + 4x(1-x)}}{2(z^2+1)}
$$
It is routine to check that $\theta_1(x,z), \theta_2(x,z) \in [0,1]$ whenever $x \in [0,1]$ and $z>0$. Whence,
\begin{align*}
&\Big{|} \int_0^1 \Phi(u(x, \theta, n)) f(\theta)\ud\theta - \ff(x)\Big{|}\\
&\quad \leq \int_0^{+\infty} \sqrt{\frac{n}{2\pi}} \exp\Big\{-\frac{n}{2}z^2\Big\} |\ff(\theta_1(x,z)) + \ff(\theta_2(x,z)) -2\ff(x)| \ud z
\end{align*}
so that, introducing $\delta_n := \sqrt{\frac{2\log (n+1)}{n}}$, we have
\begin{align}
& \int_{\delta_n}^{+\infty} \sqrt{\frac{n}{2\pi}} \exp\Big\{-\frac{n}{2}z^2\Big\} |\ff(\theta_1(x,z)) + \ff(\theta_2(x,z)) -2\ff(x)| \ud z \nonumber \\
&\leq 2 \int_{\delta_n}^{+\infty} \sqrt{\frac{n}{2\pi}} \exp\Big\{-\frac{n}{2}z^2\Big\} \ud z \leq \frac{1}{(n+1)\sqrt{\pi\log(n+1)}}\ . \label{eq:gaussext}
\end{align}
It remains to study the integral on $[0, \delta_n]$, by noticing that, after this splitting, we can consider the variable $z^2$ much smaller than $x$ and $1-x$, whenever $x \in I(n,\gamma)$. More precisely, given $\gamma \in (0,1)$ it is possible to find an integer $N(\gamma) \geq 4$ for which $\delta_n \leq \xsng(1 - \xsng)$ for all $n \geq N(\gamma)$. Therefore, we have that both 
$\overline{\theta}_1(x,z) := x - z\sqrt{x(1-x)}$ and $\overline{\theta}_2(x,z) := x + z\sqrt{x(1-x)}$ belong to $[0,1]$ whenever $z \in [0,\delta_n]$, $x \in I(n,\gamma)$ and $n \geq N(\gamma)$. 
We now have
\begin{align*}
&|\ff(\theta_1) + \ff(\theta_2) -2\ff(x)|\\
&\quad\leq |\ff(\theta_1) - \ff(\overline{\theta}_1)| + |\ff(\overline{\theta}_1) + \ff(\overline{\theta}_2) -2\ff(x)| + |\ff(\overline{\theta}_2) - \ff(\theta_2)| \nonumber \\
&\quad\leq \|f\|_{\infty}[|\theta_1 - \overline{\theta}_1| + |\theta_2 - \overline{\theta}_2|] + z^2 \Big{|} \frac{\ff(\overline{\theta}_1) - 2\ff(x) + \ff(\overline{\theta}_2)}{z^2}\Big{|} \nonumber
\end{align*}
where both $|\theta_1 - \overline{\theta}_1|$ and $|\theta_2 - \overline{\theta}_2|$ are bounded from above by $\frac{3}{2}z^2$. To check this bound, it is enough to consider the quantities
$$
\frac{|2x + z^2 \pm z\sqrt{z^2 + 4x(1-x)} - 2x(1+z^2) \mp z(1+z^2)\sqrt{4x(1-x)}|}{2(1+z^2)}
$$
which are less than $\frac{|1-2x|}{2}z^2 + \frac{1}{2}z^3 +  \frac{1}{2}z[\sqrt{z^2 + 4x(1-x)} - \sqrt{4x(1-x)}]$. The desired result now follows by observing that $|1-2x| \leq 1$, $z \in (0,1)$ whenever
$n \geq N(\gamma)$, and $\sqrt{z^2 + 4x(1-x)} - \sqrt{4x(1-x)} \leq z$. To conclude this argument, we note that $\big{|} \frac{\ff(\overline{\theta}_1) - 2\ff(x) + \ff(\overline{\theta}_2)}{z^2}\big{|}$
is bounded by virtue of \eqref{eq:der2} with $w = z\sqrt{x(1-x)}$, since $z\sqrt{x(1-x)} < x(1-x)$ holds under the restrictions $z \in [0,\delta_n]$, $x \in I(n,\gamma)$ and $n \geq N(\gamma)$.  
Whence,
\begin{align}
& \int_0^{\delta_n} \sqrt{\frac{n}{2\pi}} \exp\Big\{-\frac{n}{2}z^2\Big\} |\ff(\theta_1(x,z)) + \ff(\theta_2(x,z)) -2\ff(x)| \ud z \nonumber \\
&\quad\leq (3\|f\|_{\infty} + R(\gamma)|f|_{1,\gamma}) \int_0^{\delta_n} z^2 \sqrt{\frac{n}{2\pi}} \exp\Big\{-\frac{n}{2}z^2\Big\} \ud z \nonumber \\
&\quad= \frac{1}{2} (3\|f\|_{\infty} + R(\gamma)|f|_{1,\gamma}) \int_{\reals} z^2 \sqrt{\frac{n}{2\pi}} \exp\Big\{-\frac{n}{2}z^2\Big\} \ud z  = \frac{3\|f\|_{\infty} + R(\gamma)|f|_{1,\gamma}}{2n}\ . \ \ \ \label{eq:gaussder2} 
\end{align}
To bound the expression in \eqref{eq:splitmain} that contains $\hh_n$, we can exploit the inequalities $|1-y^2| \leq e^{y^2/3}$ and $|S(x)| \leq 1/2$, valid for any $x, y \in \reals$, 
to get 
$$
|\hh_n(y; \theta)| \leq \frac{2}{3\sqrt{2\pi n\theta(1-\theta)}} e^{-y^2/6} + \frac{\lambda_3}{n\theta(1-\theta)}
$$
for all $y \in \reals$ and $\theta \in (0,1)$, where $\lambda_3 := \frac{1}{12\sqrt{2\pi}}\sup_{y \in \reals} e^{-y^2/2}|y^3 -3y|$. Then, after writing 
\begin{equation} \label{eq:Hnfirst}
\int_0^1 |\hh_n(u(x, \theta, n); \theta)| f(\theta)\ud\theta \leq \frac{2\|f\|_{\infty}}{3\sqrt{2\pi n}} \int_0^1 \frac{e^{-u(x, \theta; n)^2/6}}{\sqrt{\theta(1-\theta)}} \ud\theta + \frac{\lambda_3}{n} \int_0^1 \frac{f(\theta)}{\theta(1-\theta)}\ud\theta\ ,
\end{equation}
we have only to show that the former integral on the right-hand side is $O(1/n)$ since, for the latter integral, it is enough to notice that 
\begin{equation} \label{eq:Hnsecond}
\int_0^1 \frac{f(\theta)}{\theta(1-\theta)}\ud\theta \leq \frac{2^{\gamma+1}}{1-\gamma}|f|_{1,\gamma} \int_0^{1/2} \theta^{-\gamma} \ud\theta + \frac{2^{\gamma+1}}{1-\gamma}|f|_{1,\gamma} \int_{1/2}^1 (1-\theta)^{-\gamma} \ud\theta \leq
\frac{2^{\gamma+2}}{(1-\gamma)^2}|f|_{1,\gamma}
\end{equation}
by virtue of the same arguments used to prove point i) of Lemma \ref{lm:regularity}. Now, the above-mentioned claim about the former integral follows after checking the boundedness of the expressions
$$
\mathfrak{I}_n(a,b;x) := \sqrt{n} \int_a^b \exp\left\{- \frac{n(x-\theta)^2}{6\theta(1-\theta)} \right\} \frac{\ud\theta}{\sqrt{\theta(1-\theta)}}
$$
by letting $n$ and $x$ vary in $\naturals$ and $I(n,\gamma)$, respectively, where $[a,b]$ coincides with either $[0, 1/2]$ or $[1/2, 1]$. Therefore, taking the former case as reference, we have
\begin{align*}
\mathfrak{I}_n\big(0,\frac12; x\big) &\leq \sqrt{\frac{n}{2}} \int_0^{\infty} \exp\left\{- \frac{n(x-\theta)^2}{6\theta} \right\} \frac{\ud\theta}{\sqrt{\theta}}\\
& = \sqrt{3} e^{nx/3} \int_0^{\infty} \exp\left\{- \frac{(\frac{1}{6}nx)^2}{y} - y\right\} \frac{\ud y}{\sqrt{y}}  \nonumber \\
&= \sqrt{2}(nx)^{1/2} e^{nx/3} K_{1/2}(nx/3) \nonumber
\end{align*}
where $K_{1/2}$ stands for the modified Bessel function of the second kind. See formula 3.471.12 in Gradshtein and Ryzik \cite{GradRyz}. Since $x \in I(n,\gamma)$ implies that $nx \geq n^{\frac{1-\gamma}{2-\gamma}} \geq 1$, we notice that
the expression $\sqrt{2}(nx)^{1/2} e^{nx/3} K_{1/2}(nx/3)$ is bounded, in view of the asymptotic expansion $K_{1/2}(z) \sim \sqrt{\frac{\pi}{2z}} e^{-z}$, which is valid as $z \rightarrow +\infty$.  
Since an analogous bound holds also for $\mathfrak{I}_n(\frac12,1;x)$, we can combine \eqref{eq:Hnfirst}-\eqref{eq:Hnsecond} with the analytical study of $\mathfrak{I}_n(a,b;x)$ to obtain that
\begin{equation} \label{eq:Hnfinal}
\sup_{x \in I(n,\gamma)} \int_0^1 |\hh_n(u(x, \theta, n); \theta)| f(\theta)\ud\theta \leq \left[ \|f\|_{\infty} + |f|_{1,\gamma} \right] \frac{\lambda_4}{n}
\ \ \ \ \ \ \ \ (n \geq 4)
\end{equation}
is valid with a numerical constant $\lambda_4$, independent of $f$ and $n$.


We conclude the first part of the proof with the analysis of the last term on the right-hand side of \eqref{eq:splitmain}. Taking account of \eqref{eq:epsn}, we immediately realize that the latter summand yields $\dfrac{\lambda_2 G}{n} 
\int_0^1 f(\theta)[\theta(1-\theta)]^{-1}\ud \theta$, which is of order $O(1/n)$ by virtue of \eqref{eq:Hnsecond}. The study of the former summand in \eqref{eq:epsn} is more laborious, and it will be conducted by mimicking the argument used in Ibragimov and Linnik \cite{IbraLinn} to prove formula (3.3.10). As first step, we borrow from Section 3.3 of Ibragimov and Linnik \cite{IbraLinn} the explicit expression of the Fourier-Stieltjes transform $d_n(t)$ (see page 101 therein) and we combine it with the formulae displayed in Section VI.1 of Petrov \cite{petrov75}, to obtain  
\begin{align}
\hat{\dd}_n(\xi; \theta)&:= \int_{\reals} e^{i\xi y} \ud_y \dd_n(y; \theta)\nonumber\\
&= \frac{-\xi}{\sqrt{n \theta(1 - \theta)}} \!\! \sum_{r \in \relatives\setminus\{0\}} \frac{e^{2\pi i r\theta n}}{2\pi r} \exp\left\{-\frac{1}{2}[\xi+2\pi r\sqrt{n\theta(1-\theta)}]^2\right\}  \nonumber \\
&\quad\times \left[1 + \frac{1-2\theta}{6\sqrt{n\theta(1 - \theta)}}(i\xi+2\pi ir\sqrt{n\theta(1-\theta)})^3 \right] \label{eq:FouStie}
\end{align}
where
$$
\dd_n(y; \theta) := \frac{1}{\sqrt{2\pi n\theta(1-\theta)}} e^{-\frac{1}{2}y^2} S(n\theta + y \sqrt{n \theta(1 - \theta)} ) \Big{[} 1 + \frac{1-2\theta}{6\sqrt{n\theta(1 - \theta)}}(y^3 - 3y) \Big{]}\ .
$$
Then, we split the integral in \eqref{eq:epsn} into five terms, by dividing the domain $[-n,n]$ into suitable subdomains whose definitions depend on $T_1(n,\theta) := \pi \sqrt{n \theta(1 - \theta)}$ and 
$$
T_2(n,\theta) := \sqrt{\frac{n \theta(1 - \theta)}{1-3\theta+3\theta^2}}\ .
$$
We observe that, since $1-3\theta+3\theta^2 \geq 1/4$ for any $\theta \in [0,1]$, the relation $T_2(n,\theta) \leq T_1(n,\theta)$ is always in force, whereas $T_1(n,\theta) \leq n$ holds whenever $n > \pi^2/4$, which we now assume. Therefore, the desired bound for the integral in \eqref{eq:epsn} follows from
\begin{align}
&\int_{-n}^n \Big{|} \frac{\hat{\bb}_n(\xi; \theta) - \hat{\gg}_n(\xi; \theta)}{\xi} \Big{|} \ud\xi \nonumber \\
&\quad\leq \int_{-T_2(n,\theta)}^{T_2(n,\theta)} \Big{|} \frac{\hat{\bb}_n(\xi; \theta) - \hat{\vv}_n(\xi; \theta)}{\xi} \Big{|} \ud\xi + 
\int_{-T_1(n,\theta)}^{T_1(n,\theta)} \Big{|} \frac{\hat{\dd}_n(\xi; \theta)}{\xi} \Big{|} \ud\xi \nonumber \\
&\quad\quad+ \int_{\{T_2(n,\theta) \leq |\xi| \leq n\}} \Big{|} \frac{\hat{\vv}_n(\xi; \theta)}{\xi} \Big{|} \ud\xi + \int_{\{T_2(n,\theta) \leq |\xi| \leq T_1(n,\theta) \}} \Big{|} \frac{\hat{\bb}_n(\xi; \theta)}{\xi} \Big{|} \ud\xi \nonumber \\
&\quad\quad+ \int_{\{T_1(n,\theta) \leq |\xi| \leq n\}} \Big{|} \frac{\hat{\bb}_n(\xi; \theta) - \hat{\dd}_n(\xi; \theta)}{\xi} \Big{|} \ud\xi \label{eq:spliteps}
\end{align}
where 
\begin{eqnarray}
\hat{\vv}_n(\xi; \theta) &:=& \int_{\reals} e^{i\xi y} \ud_y\left( \Phi(y) + \frac{1}{\sqrt{2\pi n}} e^{-\frac{1}{2}y^2} \frac{1-2\theta}{6\sqrt{\theta(1 - \theta)}}(1-y^2)\right) \nonumber \\
&=& e^{-\frac{1}{2}\xi^2} \left[1 +  \frac{1-2\theta}{6\sqrt{n\theta(1 - \theta)}}(i\xi)^3 \right] \ . \nonumber
\end{eqnarray}
For the derivation of $\hat{\vv}(\xi; \theta)$, see Section VI.1 of Petrov \cite{petrov75}. Moreover, with reference to that very same section, we note that the term $\frac{1-2\theta}{\sqrt{\theta(1 - \theta)}}$ coincides with the ratio between the third cumulant and the third power of the standard deviation of a centered Bernoulli variable with parameter $\theta$, while $y^2-1$ coincides with the Chebyshev-Hermite polynomial of degree 2. In addition, $T_2(n,\theta)$ coincides with the product between $\sqrt{n}$ and the square root of the ratio between the fourth power of the standard deviation and the fourth moment of the same centered Bernoulli variable.

In view of these remarks, we provide a bound for the first integral on the right-hand side of \eqref{eq:spliteps} by an application of Lemma 4 in Chapter VI of Petrov \cite{petrov75} with $s=4$, namely
$$
\int_{-T_2(n,\theta)}^{T_2(n,\theta)} \Big{|} \frac{\hat{\bb}_n(\xi; \theta) - \hat{\vv}_n(\xi; \theta)}{\xi} \Big{|} \ud\xi \leq \lambda_5 \frac{1-3\theta+3\theta^2}{n \theta(1 - \theta)} \int_{-T_2(n,\theta)}^{T_2(n,\theta)} (|\xi|^3+|\xi|^9) e^{-\frac{1}{12}\xi^2}\ud\xi 
$$
where $\lambda_5$ is a numerical constant specified in the proof of the quoted lemma. Whence,
\begin{equation} \label{eq:BE1}
\int_{-T_2(n,\theta)}^{T_2(n,\theta)} \Big{|} \frac{\hat{\bb}_n(\xi; \theta) - \hat{\vv}_n(\xi; \theta)}{\xi} \Big{|} \ud\xi \leq \frac{\lambda_6}{n \theta(1 - \theta)} 
\end{equation}
where $\lambda_6$ is another numerical constant (independent of $n$ and $\theta$).

Then, we study the second integral on the right-hand side of \eqref{eq:spliteps} by resorting to the explicit expression of $\hat{\dd}_n(\xi; \theta)$, to obtain
\begin{eqnarray}
\int_{-T_1(n,\theta)}^{T_1(n,\theta)} \Big{|} \frac{\hat{\dd}_n(\xi; \theta)}{\xi} \Big{|} \ud\xi &\leq& \frac{1}{\sqrt{n \theta(1 - \theta)}} \sum_{r \in \relatives\setminus\{0\}} \frac{1}{2\pi |r|} [1 + 6\pi^3 |r|^3 n\theta(1-\theta)] \times \nonumber \\
&\times& \int_{-T_1(n,\theta)}^{T_1(n,\theta)} \exp\left\{-\frac{1}{2}[\xi+2\pi r\sqrt{n\theta(1-\theta)}]^2\right\}\ud\xi \ . \label{eq:DnT1}
\end{eqnarray}
At this stage, we exploit that $r^2 - |r| \geq \frac{1}{2}r^2$ if $|r| \geq 2$ to write
\begin{eqnarray}
[\xi+2\pi r\sqrt{n\theta(1-\theta)}]^2 &\geq& \xi^2 + (2\pi r\sqrt{n\theta(1-\theta)})^2 - 4T_1(n,\theta) \pi |r| \sqrt{n\theta(1-\theta)} \nonumber \\
&=& \xi^2 + 4\pi^2 (r^2 - |r|) n\theta(1-\theta) \geq \xi^2 + 2\pi^2 r^2 n\theta(1-\theta) \ . \nonumber
\end{eqnarray}
After removing the two terms corresponding to $r = \pm 1$, the series on the right-hand side of \eqref{eq:DnT1} can be bounded by
\begin{equation} \label{eq:rgeq2}
\left(\frac{2}{\pi n\theta(1-\theta)}\right)^{1/2} \sum_{r=2}^{+\infty} \frac{1 + 6\pi^3 r^3 n\theta(1-\theta)}{r} e^{-\pi^2 n \theta(1-\theta) r^2}
\end{equation}
and then, taking cognizance that there is a suitable constant $K(\beta)$ such that $\sum_{r=2}^{+\infty} r^{\beta} e^{-\lambda r^2} \leq K(\beta) \lambda^{-(\beta+1)/2}$ holds for all $\lambda > 0$ if 
$\beta \geq 0$, we have that the expression in \eqref{eq:rgeq2} can be bounded by $\lambda_7/[n\theta(1-\theta)]$, where $\lambda_7$ is a constant (independent of $n$ and $\theta$). To handle also the terms of the series corresponding to 
$r = \pm 1$, we take account that $(x \pm 2)^2 \geq \frac{1}{2}x^2 + \frac{1}{3}$ holds for all $x \in [-1,1]$, to write
$[\xi \pm 2\pi \sqrt{n\theta(1-\theta)}]^2 \geq \frac{1}{2}\xi^2 + \frac{\pi^2}{3}n\theta(1 - \theta)$ for all $\xi \in [-T_1(n,\theta), T_1(n,\theta)]$. Lastly, the sum of the two terms in \eqref{eq:DnT1} corresponding to $r = \pm 1$ can be bounded by
$$
2\left(\frac{1}{\pi n\theta(1-\theta)}\right)^{1/2} [1 + 6\pi^3 n\theta(1-\theta)] e^{-\frac{\pi^2}{6} n \theta(1-\theta)}\ .
$$
Therefore, we can conclude that 
\begin{equation} \label{eq:BE2}
\int_{-T_1(n,\theta)}^{T_1(n,\theta)} \Big{|} \frac{\hat{\dd}_n(\xi; \theta)}{\xi} \Big{|} \ud\xi \leq \frac{\lambda_8}{n \theta(1 - \theta)} 
\end{equation}
holds with a suitable numerical constant $\lambda_8$ (independent of $n$ and $\theta$). 

As for the third integral on the right-hand side of \eqref{eq:spliteps}, we just use the explicit expression of $\hat{\vv}_n(\xi; \theta)$ to write  
\begin{align*}
&\int_{\{T_2(n,\theta) \leq |\xi| \leq n\}} \Big{|} \frac{\hat{\vv}_n(\xi; \theta)}{\xi} \Big{|} \ud\xi\\
&\quad\leq 2 \int_{\sqrt{n\theta(1-\theta)}}^{+\infty} \frac{e^{-\frac{1}{2}\xi^2}}{\xi} \ud\xi + \frac{1}{3\sqrt{n\theta(1 - \theta)}}
\!\!\!\int_{\sqrt{n\theta(1-\theta)}}^{+\infty} \xi^2 e^{-\frac{1}{2}\xi^2} \ud\xi \ .
\end{align*}
Using that $x^p e^{-x} \leq (p/e)^p$, which is valid whenever $x,p > 0$, we show that  
\begin{equation} \label{eq:BE3}
\int_{\{T_2(n,\theta) \leq |\xi| \leq n\}} \Big{|} \frac{\hat{\vv}_n(\xi; \theta)}{\xi} \Big{|} \ud\xi \leq \frac{\lambda_9}{n \theta(1 - \theta)} 
\end{equation}
holds with a suitable numerical constant $\lambda_9$ (independent of $n$ and $\theta$). 

We now consider the fourth integral on the right-hand side of \eqref{eq:spliteps}. By definition, we have
\begin{gather}
\int_{\{T_2(n,\theta) \leq |\xi| \leq T_1(n,\theta) \}} \Big{|} \frac{\hat{\bb}_n(\xi; \theta)}{\xi} \Big{|} \ud\xi \nonumber \\
= \int_{\{T_2(n,\theta) \leq |\xi| \leq T_1(n,\theta) \}} |\xi|^{-1}\ \Big{|} \ee\Big{[} \exp\Big{\{} \frac{i\xi}{\sqrt{n\theta(1-\theta)}} \Big{(} \sum_{j=1}^n X_j - n\theta \Big{)}  \Big{\}} \ |\ Y=\theta \Big{]} \Big{|}  \ud\xi \nonumber
\end{gather}
and, after changing the variable by the rule $u = \xi/\sqrt{n\theta(1-\theta)}$ and noticing that $(1-3\theta+3\theta^2)^{-1/2} \geq 1$ is valid for any $\theta \in [0,1]$, we provide the following upper bound
$$
\int_{\{1 \leq |u| \leq \pi \}} \Big{|} \ee\Big{[} \exp\{ iu\sum_{j=1}^n X_j \} \ |\ Y=\theta \Big{]} \Big{|}  \ud u\ .
$$
Now, we just utilize the explicit form of the characteristic function of the binomial distribution with parameters $n$ and $\theta$ to write
$$
\Big{|} \ee\Big{[} \exp\{ iu\sum_{j=1}^n X_j \} \ |\ Y=\theta \Big{]} \Big{|} = |1 - \theta + \theta e^{iu}|^n = [1 - 2\theta(1-\theta)(1-\cos u)]^{n/2}\ .
$$
Whence,
\begin{eqnarray} 
\int_{\{T_2(n,\theta) \leq |\xi| \leq T_1(n,\theta) \}} \Big{|} \frac{\hat{\bb}_n(\xi; \theta)}{\xi} \Big{|} \ud\xi &\leq& 2(\pi-1) [1 - 2\theta(1-\theta)(1-\cos 1)]^{n/2} \nonumber \\
&\leq& \frac{2(\pi-1)}{e(1-\cos 1)} \frac{1}{n\theta(1-\theta)}\ . \label{eq:BE4}
\end{eqnarray}

To study of the last integral on the right-hand side of \eqref{eq:spliteps}, we introduce the characteristic function $\phi(\cdot; \theta)$ of the r.v. $(X_1 - \theta)$ given $Y=\theta$, that is
$\phi(\xi; \theta) = [(1-\theta) + \theta e^{i\xi}] e^{-i\xi\theta}$, so that we have $\hat{\bb}_n(\xi; \theta) = [\phi(\xi/\sqrt{n\theta(1-\theta)}; \theta)]^n$. After changing the variable in that integral according to 
$u = \xi/\sqrt{n\theta(1-\theta)}$ and recalling that $\phi(-\xi; \theta) = \overline{\phi(\xi; \theta)}$, we get
\begin{align*}
 &\int_{\{T_1(n,\theta) \leq |\xi| \leq n\}} \Big{|} \frac{\hat{\bb}_n(\xi; \theta) - \hat{\dd}_n(\xi; \theta)}{\xi} \Big{|} \ud\xi\\
 &\quad= 2 \int_{\pi}^{\sqrt{n/[\theta(1-\theta)]}}
 \Big{|} \frac{[\phi(u; \theta)]^n - \hat{\dd}_n(u\sqrt{n\theta(1-\theta)}; \theta)}{u} \Big{|} \ud u\ .
\end{align*}
At this stage, we introduce the quantity
$$
\overline{r}(n;\theta) := \left\lfloor \frac{1}{2}\left(\frac{1}{\pi}\sqrt{\frac{n}{\theta(1-\theta)}} -1 \right) \right\rfloor 
$$
and we notice that $(2\overline{r}(n;\theta) + 1)\pi \leq \sqrt{\frac{n}{\theta(1-\theta)}} < (2\overline{r}(n;\theta) + 3)\pi$. In this notation, we have
\begin{equation} \label{eq:IL3313}
\int_{\pi}^{\sqrt{n/[\theta(1-\theta)]}} \Big{|} \frac{[\phi(u; \theta)]^n - \hat{\dd}_n(u\sqrt{n\theta(1-\theta)}; \theta)}{u} \Big{|} \ud u
\leq \sum_{k=1}^{\overline{r}(n;\theta) + 1} \mathfrak{J}_k(n;\theta)
\end{equation}
where
$$
\mathfrak{J}_k(n;\theta) := \int_{(2k-1)\pi}^{(2k+1)\pi} \Big{|} \frac{[\phi(u; \theta)]^n - \hat{\dd}_n(u\sqrt{n\theta(1-\theta)}; \theta)}{u} \Big{|} \ud u\ .
$$
To bound the integrals $\mathfrak{J}_k$'s, we first isolate from the series \eqref{eq:FouStie} defining $\hat{\dd}_n(u\sqrt{n\theta(1-\theta)}; \theta)$ the term corresponding to $r = -k$, which reads
\begin{align*}
&\frac{u e^{-2\pi i kn\theta}}{2k\pi} \exp\left\{-\frac{n\theta(1-\theta)}{2}(u-2k\pi)^2\right\} \cdot \left[1 + \frac{1-2\theta}{6}n\theta(1-\theta)(iu - 2\pi k i)^3 \right]\\
&\quad =: \hat{\delta}_{n,k}(u;\theta),
\end{align*}
so that we obtain
\begin{align*}
\mathfrak{J}_k(n;\theta) &\leq \int_{(2k-1)\pi}^{(2k+1)\pi} \Big{|} \frac{\hat{\dd}_n(u\sqrt{n\theta(1-\theta)}; \theta) - \hat{\delta}_{n,k}(u;\theta)}{u} \Big{|} \ud u\\
&\quad + \int_{(2k-1)\pi}^{(2k+1)\pi} \Big{|} \frac{[\phi(u; \theta)]^n - \hat{\delta}_{n,k}(u;\theta)}{u} \Big{|} \ud u \nonumber \\ 
&\quad\quad=: \mathfrak{J}_k^{(1)}(n;\theta)  + \mathfrak{J}_k^{(2)}(n;\theta) \nonumber \ . 
\end{align*}
To analyze $\mathfrak{J}_k^{(1)}(n;\theta)$, we write
\begin{align} 
\mathfrak{J}_k^{(1)}(n;\theta) &\leq \left(\sum_{r=1}^{\infty} + \sum_{r=-\infty}^{-(k+2)} + \sum_{r\in\{-(k+1)\}} +  \sum_{r\in\{-k+1\}} + \sum_{r =-k+2}^{-1} \right)  \frac{1}{2\pi |r|} \nonumber \\
&\quad\times \int_{(2k-1)\pi}^{(2k+1)\pi} \exp\left\{-\frac{n\theta(1-\theta)}{2} (u+2\pi r)^2 \right\}\nonumber \\
&\quad\quad\times \left[1 + \frac{1}{6}n\theta(1 - \theta) |u + 2\pi r|^3\right] \ud u \label{eq:J1} \ \ \ \ \ 
\end{align}
with the proviso that both the fourth and the fifth sum are void when $k=1$, and that the fifth sum is void when $k=2$.

To deal with the series in \eqref{eq:J1} limited to $r \in \naturals$, we observe that $(u+2\pi r)^2 \geq u^2 + 4\pi^2r^2$ if $u \in [(2k-1)\pi, (2k+1)\pi]$ and we take account that $|z_1 + z_2|^3 \leq 4(|z_1|^3 + |z_2|^3)$, so that we deduce,
for the series at issue, the upper bound
\begin{align*}
&\sum_{r=1}^{\infty} \frac{1}{2\pi r} e^{-2\pi^2n\theta(1-\theta)r^2}\\
&\times  \int_{(2k-1)\pi}^{(2k+1)\pi} \exp\left\{-\frac{n\theta(1-\theta)}{2} u^2 \right\} \left[1 + \frac{2}{3}n\theta(1 - \theta)(u^3 + 8\pi^3r^3)\right] \ud u\ .
\end{align*}
At this stage, recalling \eqref{eq:IL3313}, we conclude that the sum over the index $k$ of the last expression is majorized by 
\begin{align*}
& \sum_{r=1}^{\infty} \frac{1}{2\pi r} e^{-2\pi^2n\theta(1-\theta)r^2}  \int_0^{+\infty} \exp\left\{-\frac{n\theta(1-\theta)}{2} u^2 \right\} \left[1 + \frac{2}{3}n\theta(1 - \theta)(u^3 + 8\pi^3r^3)\right] \ud u \nonumber \\
&\quad= \sum_{r=1}^{\infty} \frac{1}{2\pi r} e^{-2\pi^2n\theta(1-\theta)r^2}  \left[ \frac{1}{2}\sqrt{\frac{2\pi}{n\theta(1 - \theta)}} + \frac{4}{3n\theta(1-\theta)} +  \frac{8\sqrt{2\pi}}{3}\pi^3r^3\sqrt{n\theta(1-\theta)} \right] \ . \nonumber
\end{align*}
Then, we use $x^p e^{-x} \leq (p/e)^p$, valid for any $x, p > 0$, with $p = \frac{1}{2} + \varepsilon(\gamma), \varepsilon(\gamma), \frac{3}{2} + \varepsilon(\gamma)$, respectively, and $\varepsilon(\gamma) := \frac{1-\gamma}{2}$,
to produce the global bound $S_1(\gamma) [n\theta(1-\theta)]^{-(1+\varepsilon(\gamma))}$ for the last series, where
\begin{eqnarray}
S_1(\gamma) &:=& \frac{1}{2\pi}  \Big[ \frac{\sqrt{2\pi}}{2} \left(\frac{\frac{1}{2} + \varepsilon(\gamma)}{2\pi^2 e}\right)^{\frac{1}{2} + \varepsilon(\gamma)} \!\!\!\!\!\! \zeta(2[1+\varepsilon(\gamma)]) + 
\frac{4}{3} \left(\frac{\varepsilon(\gamma)}{2\pi^2 e}\right)^{\varepsilon(\gamma)} \!\!\!\!\! \zeta(1+2\varepsilon(\gamma)) \nonumber \\
&+&  \frac{8\pi^3\sqrt{2\pi}}{3} \left(\frac{\frac{3}{2} + \varepsilon(\gamma)}{2\pi^2 e}\right)^{\frac{3}{2} + \varepsilon(\gamma)} \!\!\!\!\!\! \zeta(1+2\varepsilon(\gamma)) \Big] \nonumber 
\end{eqnarray}
$\zeta(\cdot)$ denoting the Riemann zeta function. Now, we come back to \eqref{eq:J1} and we consider the remaining four sums. The change of variable $u = s + 2k\pi$ in the integral and the inequality $|z_1 + z_2|^3 \leq 4(|z_1|^3 + |z_2|^3)$ lead us to rewrite the expression inside the sums in \eqref{eq:J1} as
\begin{equation} \label{eq:moldava}
\frac{1}{2\pi |r|} [1 + 6\pi^3n\theta(1 - \theta)|r+k|^3]  \int_{-\pi}^{\pi} \exp\left\{-\frac{n\theta(1-\theta)}{2} [s+2\pi(k+r)]^2 \right\} \ud s \ .
\end{equation}
Therefore, for the second series in \eqref{eq:J1}, relative to the set $r \leq -(k+2)$, we have
\begin{align*}
&\sum_{r=k+2}^{\infty} \frac{1}{2\pi r} [1 + 6\pi^3n\theta(1 - \theta)|k-r|^3] \int_{-\pi}^{\pi} \exp\left\{-\frac{n\theta(1-\theta)}{2} [s+2\pi(k-r)]^2 \right\} \ud s \nonumber \\
&\quad\leq \frac{1}{2\pi k} \sum_{h=2}^{\infty} \left[1 + 6\pi^3n\theta(1 - \theta)h^3 \right] \int_{-\pi}^{\pi} \exp\left\{-\frac{n\theta(1-\theta)}{2} (s-2\pi h)^2 \right\} \ud s\ . \nonumber
\end{align*}
After noticing that $(s-2\pi h)^2 \geq s^2 + 2\pi^2h^2$ for $s \in [-\pi, \pi]$ and $h \geq 2$, we get the new upper bound
$$
\frac{1}{2\pi k}  \int_{-\pi}^{\pi} \exp\left\{-\frac{n\theta(1-\theta)}{2} s^2 \right\} \ud s \times \sum_{h=2}^{\infty} \left[1 + 6\pi^3n\theta(1 - \theta)h^3 \right] e^{- \pi^2 n\theta(1-\theta) h^2}
$$
which is less or equal than $\frac{1}{k} S_2(\gamma) [n\theta(1-\theta)]^{-(1+\varepsilon(\gamma))}$ where, by another application of $x^p e^{-x} \leq (p/e)^p$ for $p = \frac{1}{2} + \varepsilon(\gamma)$ and $2 + \varepsilon(\gamma)$, respectively, 
$$
S_2(\gamma) := \frac{1}{2\pi} \Big[ \sqrt{2\pi} \left(\frac{\frac{1}{2} + \varepsilon(\gamma)}{\pi^2 e}\right)^{\frac{1}{2} + \varepsilon(\gamma)} \!\!\!\!\!\! \zeta(1+2\varepsilon(\gamma)) + 12\pi^4 \left(\frac{2 + \varepsilon(\gamma)}{\pi^2 e}\right)^{2 
+ \varepsilon(\gamma)} \!\!\!\!\!\! \zeta(2[1+\varepsilon(\gamma)]) \Big] \ . 
$$
For $r = -(k+1)$ the expression in \eqref{eq:moldava} is majorized by
$$
\frac{1 + 6\pi^3n\theta(1 - \theta)}{k+1} \exp\left\{-\frac{\pi^2 n\theta(1-\theta)}{2} \right\}  
$$
which is, in turn, less or equal than $\frac{1}{k+1} S_3(\gamma) [n\theta(1-\theta)]^{-(1+\varepsilon(\gamma))}$ with
$$
S_3(\gamma) :=  \left(\frac{2[1 + \varepsilon(\gamma)]}{\pi^2 e}\right)^{1 + \varepsilon(\gamma)} +\  6\pi^3 \left(\frac{2[2 + \varepsilon(\gamma)]}{\pi^2 e}\right)^{2 + \varepsilon(\gamma)}\ .
$$
Analogously, for any $k \geq 2$, the expression in \eqref{eq:moldava} with $r = -k+1$ is majorized by
$$
\frac{1 + 6\pi^3n\theta(1 - \theta)}{k-1} \exp\left\{-\frac{\pi^2 n\theta(1-\theta)}{2} \right\}  
$$
which is, in turn, less or equal than $\frac{1}{k-1} S_3(\gamma) [n\theta(1-\theta)]^{-(1+\varepsilon(\gamma))}$. Finally, for $k \geq 3$, it remains to provide an upper bound for the sum of the expression \eqref{eq:moldava} as $r$ varies from
$-k+2$ to $-1$. Changing the variable in the sum, according to $h=k+r$, we obtain the equivalent expression
$$
\sum_{h=2}^{k-1} \frac{1}{2\pi (k-h)} [1 + 6\pi^3n\theta(1 - \theta)h^3]  \int_{-\pi}^{\pi} \exp\left\{-\frac{n\theta(1-\theta)}{2} (s+2\pi h)^2 \right\} \ud s 
$$
which is majorized by virtue of the inequality $(s+2\pi h)^2 \geq s^2 + 2\pi^2h^2$, valid for any $s \in [-\pi,\pi]$ and $h \geq 2$. Then, we arrive at
$$
\int_{-\pi}^{\pi} \exp\left\{-\frac{n\theta(1-\theta)}{2} s^2 \right\} \ud s \times \sum_{h=2}^{k-1} \frac{e^{-\pi^2n\theta(1-\theta)h^2}}{2\pi (k-h)} [1 + 6\pi^3n\theta(1 - \theta)h^3] 
$$
and we now realize that, in view of \eqref{eq:IL3313}, we can exchange the order of summation according to $\sum_{k=3}^{\overline{r}(n;\theta) + 1} \sum_{h=2}^{k-1} = \sum_{h=2}^{\overline{r}(n;\theta)} \sum_{k=h+1}^{\overline{r}(n;\theta) + 1}$.
At this stage, for the inner sum, we have
$$
\sum_{k=h+1}^{\overline{r}(n;\theta) + 1} \frac{1}{k-h} \leq \sum_{r=1}^{\overline{r}(n;\theta)} \frac{1}{r} \leq \frac{1}{2} \log\left(\frac{n}{\theta(1-\theta)} \right)\ .
$$
Using this upper bound, we pass to the outer sum, which is majorized by 
\begin{align*}
&\int_{-\pi}^{\pi} \exp\left\{-\frac{n\theta(1-\theta)}{2} s^2 \right\} \ud s\\
&\quad \times \frac{1}{4\pi} \log\left(\frac{n}{\theta(1-\theta)} \right) \sum_{h=2}^{\infty} e^{-\pi^2n\theta(1-\theta)h^2} [1 + 6\pi^3n\theta(1 - \theta)h^3] \ .
\end{align*}
The series in the above expression has been already treated above, yielding the further upper bound
$$
\frac{S_2(\gamma)}{2} \log\left(\frac{n}{\theta(1-\theta)} \right) [n\theta(1-\theta)]^{-(1+\varepsilon(\gamma))}\ .
$$
Therefore, gathering all the bounds that follow formula \eqref{eq:J1}, we get
\begin{equation} \label{eq:J1final}
\sum_{k=1}^{\overline{r}(n;\theta) + 1} \mathfrak{J}_k^{(1)}(n;\theta) \leq \frac{S_1(\gamma)}{ [n\theta(1-\theta)]^{(1+\varepsilon(\gamma))}} + \frac{S_2(\gamma)+S_3(\gamma)} { [n\theta(1-\theta)]^{(1+\varepsilon(\gamma))}}
\log\left(\frac{n}{\theta(1-\theta)} \right)\ .
\end{equation}

Now, we pass to analyze the integrals $\mathfrak{J}_k^{(2)}$'s. We start again from the change of variable $u = s + 2k\pi$ and we exploit the fact that $\phi(s + 2k\pi; \theta) = \phi(s; \theta) e^{-2\pi i k\theta}$, to obtain
\begin{align*}
& \int_{(2k-1)\pi}^{(2k+1)\pi} \big{|} \frac{[\phi(u; \theta)]^n - \hat{\delta}_{n,k}(u;\theta)}{u} \big{|} \ud u \nonumber \\
&\quad= \int_{-\pi}^{\pi} \frac{\big{|}[\phi(s; \theta)]^n - \left(1 + \frac{s}{2k\pi}\right) \exp\left\{-\frac{n\theta(1-\theta)}{2}s^2 \right\} \left[1 + \frac{1-2\theta}{6}n\theta(1-\theta)(is)^3 \right]  \big{|}}{s + 2k\pi} \ud s\ . \nonumber
\end{align*}
The last integral is majorized by
\begin{align}
& \frac{1}{\pi (2k-1)}\int_{-\pi}^{\pi} \Big{|}[\phi(s; \theta)]^n - \exp\left\{-\frac{n\theta(1-\theta)}{2}s^2 \right\} \left[1 + \frac{1-2\theta}{6}n\theta(1-\theta)(is)^3 \right] \Big{|}\ \ud s \nonumber \\
&+ \frac{1}{\pi^2 k(2k-1)} \int_0^{\infty} s \exp\left\{-\frac{n\theta(1-\theta)}{2}s^2 \right\} \left[1 + \frac{\pi}{6}n\theta(1-\theta)s^2 \right] \ud s\ .
\label{eq:IL3314}
\end{align}
For the first summand in \eqref{eq:IL3314}, we change again the variable according to $s = \xi/\sqrt{n\theta(1-\theta)}$ to obtain the equality with 
$$
\frac{1}{\pi (2k-1)\sqrt{n\theta(1-\theta)}} \int_{-T_1(n,\theta)}^{T_1(n,\theta)} \big{|} \hat{\bb}_n(\xi; \theta) - \hat{\vv}_n(\xi; \theta) \big{|} \ud\xi 
$$
which can be bounded, as before, the splitting the above integral as
\begin{align*}
&\int_{-T_2(n,\theta)/4}^{T_2(n,\theta)/4} \big{|} \hat{\bb}_n(\xi; \theta) - \hat{\vv}_n(\xi; \theta) \big{|} \ud\xi\\\ 
&\quad+\ 2\int_{T_2(n,\theta)/4}^{T_1(n,\theta)} \big{|} \hat{\vv}_n(\xi; \theta) \big{|} \ud\xi \\
&\quad\quad+\ 2\int_{T_2(n,\theta)/4}^{T_1(n,\theta)} \big{|} \hat{\bb}_n(\xi; \theta) \big{|} \ud\xi \ .
\end{align*}
In fact, it is now crucial to observe that the expression 
$$
T_3(n,\theta,\delta) := \frac{1}{4}\sqrt{n} \left( \frac{[\theta(1-\theta)]^{(3+\delta)}/2}{\theta(1-\theta)[(1-\theta)^{2+\delta}+\theta^{(2+\delta)}]} \right)^{1/(1+\delta)}\ ,
$$
corresponding to the limitation for $|\xi|$ given in Lemma \ref{lm:BE} when the $V_n$'s are i.i.d., centered Bernoulli variables, is not less than $T_2(n,\theta)/4$, by virtue of the H\"older inequality. In fact, it is enough to observe that, for a centered r.v. $V$, we have $\sqrt{\sigma^4/\beta_4} \leq (\sigma^{3+\delta}/\beta_{3+\delta})^{1/(1+\delta)}$ where $\sigma^2 := \ee[V^2]$ and $\beta_s := \ee[|V|^s]$. Therefore, we can apply Lemma \ref{lm:BE}
with $\delta = 2\varepsilon(\gamma)=1-\gamma$, to get
$$
\int_{-T_2(n,\theta)/4}^{T_2(n,\theta)/4} \big{|} \hat{\bb}_n(\xi; \theta) - \hat{\vv}_n(\xi; \theta) \big{|} \ud\xi \leq \frac{\lambda_{10}}{[n\theta(1-\theta)]^{1/2 + \varepsilon(\gamma)}}\ .
$$
Since an analogous bound is in force also for $\int_{T_2(n,\theta)/4}^{T_1(n,\theta)} \big{|} \hat{\vv}_n(\xi; \theta) \big{|} \ud\xi$ and for $\int_{T_2(n,\theta)/4}^{T_1(n,\theta)} \big{|} \hat{\bb}_n(\xi; \theta) \big{|} \ud\xi$, in view of the argument already used to prove \eqref{eq:BE3}-\eqref{eq:BE4}, we conclude that
\begin{gather}
\sum_{k=1}^{\overline{r}(n;\theta) + 1} \frac{1}{\pi (2k-1)}\int_{-\pi}^{\pi} \Big{|}[\phi(s; \theta)]^n - \exp\left\{-\frac{n\theta(1-\theta)}{2}s^2 \right\} \cdot \left[1 + \frac{1-2\theta}{6}n\theta(1-\theta)(is)^3 \right] \Big{|}\ \ud s \nonumber \\
\leq \frac{\lambda_{11}}{[n\theta(1-\theta)]^{1+\varepsilon(\gamma)}} \log\left(\frac{n}{\theta(1-\theta)}\right)\ . \nonumber
\end{gather}
As to the latter summand in \eqref{eq:IL3314}, it is enough to notice that it equals
$$
\frac{1}{\pi^2 k(2k-1)} \left(2 + \frac{\pi}{12}\right) \frac{1}{n\theta(1-\theta)}
$$
yielding in the end that
\begin{equation} \label{eq:J2final}
\sum_{k=1}^{\overline{r}(n;\theta) + 1} \mathfrak{J}_k^{(1)}(n;\theta) \leq  \frac{\lambda_{12}}{n\theta(1-\theta)} + \frac{\lambda_{13}}{[n\theta(1-\theta)]^{1+\varepsilon(\gamma)}} \log\left(\frac{n}{\theta(1-\theta)}\right)\ .
\end{equation}

At this stage, we notice that, for any $\eta > 0$, $[\theta(1-\theta)]^{\eta} \log\left(\frac{1}{\theta(1-\theta)}\right)$ is bounded by a constant which depends only on $\eta$, and we can choose $\eta = \eta(\gamma) = (1-\gamma)/4$. Then, 
we collect \eqref{eq:spliteps}-\eqref{eq:BE1}-\eqref{eq:BE2}-\eqref{eq:BE3}-\eqref{eq:BE4}-\eqref{eq:IL3313}-\eqref{eq:J1final}-\eqref{eq:J2final} to draw the important conclusion that
\begin{equation} \label{eq:epsfinal}
\int_0^1 \epsilon_n(\theta) f(\theta)\ud\theta \leq \lambda_{14} \frac{\|f\|_{\infty} + |f|_{1,\gamma}}{n}
\end{equation}
holds with a suitable constant $\lambda_{14}$ which is independent of $n$ and $f$, thanks to point i) in Lemma \ref{lm:regularity} and
$$
\int_0^1 \frac{1}{[\theta(1-\theta)]^{1+\varepsilon(\gamma)+\eta(\gamma)}} f(\theta) \ud\theta \leq \frac{2^{3+\gamma+\varepsilon(\gamma)+\eta(\gamma)}}{1-\gamma} M(f)\ . 
$$ 
Therefore, the achievement of the bound \eqref{eq:epsfinal} concludes the first part of the proof, culminating in the validity of \eqref{eq:MAIN} with a suitable constant $C(\mu)$ proportional to $1 + \|f\|_{\infty} + |f|_{1,\gamma}$, under the 
the additional hypothesis $f(0)=f(1)=0$, thanks to the combination of \eqref{eq:splitfirst}-\eqref{eq:extremeff}, the two bound $9M(f)/n$ for both $\ff_n(\xsng)$ and $1-\ff_n(1-\xsng)$, and \eqref{eq:splitmain}-\eqref{eq:gaussext}-\eqref{eq:gaussder2}-\eqref{eq:Hnfinal}-\eqref{eq:epsfinal}. For completeness, we note that we have proved \eqref{eq:MAIN} only for $n\geq N_{\ast}:=\max\{4, n_0, N(\gamma), \lfloor\pi^2/4\rfloor+1\}$, but now it is immediate to extend the validity of \eqref{eq:MAIN} to all the set of positive integer: we just add the term $N_{\ast}/n$ to the right-hand side of \eqref{eq:MAIN} and we rename the new constant as $C(\mu)$.

After proving the theorem under the additional hypothesis $f(0)=f(1)=0$, we show how to get rid of this extra-condition. First, we assume that $f$ is given by a polynomial, with generic values of $f(0)$ and $f(1)$, and we apply Lemma \ref{lm:decomposition}. Since $\ff_n(x) =  \pp[S_n \leq nx] =  \int^1_0 \pp[S_n \leq nx\ |\ Y = \theta] \mu(\ud\theta)$, we obtain $\ff_n(x) = A_{\infty}\ff_{\infty,n}(x) + A_{+}\ff_{+,n}(x) - A_{-}\ff_{-,n}(x)$ for all $x \in [0,1]$, where $\ff_{\star,n}(x) := \int^1_0 \pp[S_n \leq nx\ |\ Y = \theta] f_{\star}(\theta) \ud\theta$, for $\star = \infty, +$ and $-$, respectively. Whence,
\begin{eqnarray}
\ud_K(\mu_n; \mu) &\leq&\ A_{\infty} \sup_{x \in [0,1]} |\ff_{\infty,n}(x) - \ff_{\infty}(x)|\ +\ A_{+} \sup_{x \in [0,1]} |\ff_{+,n}(x) - \ff_{+}(x)| \nonumber \\
&+&\ A_{-} \sup_{x \in [0,1]} |\ff_{-,n}(x) - \ff_{-}(x)| \label{eq:decomposition}
\end{eqnarray}
where $\ff_{\star}(x) := \int_0^x f_{\star}(\theta) \ud\theta$ for all $x \in [0,1]$ and $\star = \infty, +$ and $-$, respectively. At this stage, from Theorem 2 in Mnatsakanov \cite{Mna},
we can find a constant $C(f_{\infty})$, proportional to $1 + \|f_{\infty}\|_{\infty} + \|f'_{\infty}\|_{\infty}$, such that $\sup_{x \in [0,1]} |\ff_{\infty,n}(x) - \ff_{\infty}(x)| \leq C(f_{\infty})/n$ is in force for all $n \in \naturals$.
For the last two terms on the right-hand side of \eqref{eq:decomposition}, since $f_{\pm}(0) = f_{\pm}(1) = 0$, the problem is traced back to the first part of the proof. Hence, the inequality
\eqref{eq:MAIN} holds for all $n \in \naturals$ with a constant $C(\mu)$ proportional to $1 + \|f\|_{\infty} + |f|_{1,\gamma}$, thanks to the bounds provided in points i)-ii)-iii) of Lemma \ref{lm:decomposition}.
 
 The final act consists in removing the regularity of $f$ by some approximation arguments. First, we start from a probability density $f$ belonging to $\mathrm{C}^1([0,1])$ and we consider an approximating family of probability densities 
$f^{(\delta)}$ expressed by a polynomial which converges to $f$ uniformly with the first derivative, as $\delta \rightarrow 0$. See, e.g., Lorentz \cite{lorentz} for classical results about this kind of approximation. Since $\|f^{(\delta)}\|_{\infty} \rightarrow
\|f\|_{\infty}$ and $|f^{(\delta)}|_{1,\gamma} \rightarrow |f|_{1,\gamma}$ are obvious, we pass to analyze the behavior of $\ud_K(\mu_n; \mu)$ under the approximation. After fixing $n$, for any $x \not\in\{0, \frac{1}{n}, \dots, 1\}$, we have
\begin{align*}
|\ff_n(x) - \ff(x)| &= \lim_{\delta \rightarrow 0} |\ff_n^{(\delta)}(x) - \ff^{(\delta)}(x)|\\
&\leq \lim_{\delta \rightarrow 0}  C\frac{1 + \|f^{(\delta)}\|_{\infty} + |f^{(\delta)}|_{1,\gamma}}{n} = C\frac{1 + \|f\|_{\infty} + |f|_{1,\gamma}}{n},
\end{align*}
where the inequality follows from the previous argument. This relation entails the validity of \eqref{eq:MAIN} for all $n \in \naturals$ and $f \in \mathrm{C}^1([0,1])$, with a constant $C(\mu)$ proportional to $1 + \|f\|_{\infty} + |f|_{1,\gamma}$.
Finally, the removal of the $\mathrm{C}^1([0,1])$-regularity follows by standard arguments based on the convolution of a regularizing kernel.


\appendix
\section*{Appendix}
\numberwithin{equation}{section}
\numberwithin{thm}{section}
\numberwithin{lem}{section}
\numberwithin{prp}{section}

\subsection{Proof of Proposition \ref{prop:beta}}

We start by dealing with the case of a beta distribution with parameters $(\alpha, 1)$. First, we note that the associated density belongs to $\mathrm{W}^{1,\infty}(0,1)$ if $\alpha \in \{1\} \cup [2, +\infty)$, so that \eqref{eq:beta} follows as a direct application of Theorem 2 in Mnatsakanov \cite{Mna}. Therefore, we treat the case $\alpha \in (0,1) \cup (1,2)$ by starting from the direct computation of the d.f. $\ff_n$ associated to $\mu_n$, namely
\begin{align} 
\ff_n(x) &= \sum_{k=0}^{\lfloor nx \rfloor} \binom{n}{k} \int_0^1 \theta^k (1-\theta)^{n-k} \mu(\ud\theta) = \int^1_0 \beta(y; \lfloor nx \rfloor + 1, n - \lfloor nx \rfloor) \ff(y) \ud y \label{eq:feller} \\
&= \frac{\Gamma(n+1)}{\Gamma(n+1+\alpha)} \frac{\Gamma(\lfloor nx \rfloor + 1+\alpha)}{\Gamma(\lfloor nx \rfloor+1)} \label{eq:fellergamma}
\end{align}
for all $x \in (0,1)$, where $\beta$ is the same as in \eqref{beta_density} and 
$\lfloor\cdot\rfloor$ denotes the integral part. For the validity of the second identity in \eqref{eq:feller}, see formulae (13)-(14) in Mnatsakanov \cite{Mna}, or Problems 44-45 at the end of Chapter VI of Feller \cite{feller}. Now, for $\alpha \in (0,1)$, we invoke \emph{Wendell's inequalities} (see formula (5) in Qi and Luo \cite{Qi}) to obtain
$$
L_{\alpha}(\lfloor nx\rfloor,n) \leq \ff_n(x) \leq U_{\alpha}(\lfloor nx\rfloor,n)
$$
for all $n \in \naturals$ and $x \in (0,1)$, where $L_{\alpha}(\lfloor nx\rfloor,n) := \left(\frac{\lfloor nx\rfloor+1}{\lfloor nx\rfloor+1+\alpha}\right)^{1-\alpha} \left(\frac{\lfloor nx\rfloor+1}{n+1}\right)^{\alpha}$ and $U_{\alpha}(\lfloor nx\rfloor,n):= \left(\frac{n+1+\alpha}{n+1}\right)^{1-\alpha}\left(\frac{\lfloor nx\rfloor+1}{n+1}\right)^{\alpha}$. Then, we observe that we can write
\begin{align*}
\ud_K(\mu_n; \mu) &= \max_{k \in \{1, \dots, n\}} \sup_{x \in [\frac{k-1}{n}, \frac{k}{n})} |\ff_n(x) - \ff(x)|\\
& = \max_{k \in \{1, \dots, n\}} \sup_{x \in [\frac{k-1}{n}, \frac{k}{n})} \Big{|} \ff_n\left(\frac{k-1}{n}\right) - x^{\alpha} \Big{|} 
\end{align*}
and that, for any $k \in \{1, \dots, n\}$, 
\begin{align}
&\sup_{x \in [\frac{k-1}{n}, \frac{k}{n})}\!\! \Big{|} \ff_n\left(\frac{k-1}{n}\right) - x^{\alpha} \Big{|}\!\!  \nonumber \\
&\quad\leq [U_{\alpha}(k-1,n) - L_{\alpha}(k-1,n)] + \Big{|} U_{\alpha}(k-1,n) - \left(\frac{k}{n+1}\right)^{\alpha}\! \Big{|} \nonumber \\
&\quad\quad+ \left[ \left(\frac{k}{n}\right)^{\alpha} - \left(\frac{k}{n+1}\right)^{\alpha} \right] + \left[ \left(\frac{k}{n}\right)^{\alpha} - \left(\frac{k-1}{n}\right)^{\alpha} \right]\ . \label{eq:UnLn} 
\end{align}
For the first two summands on the above right-hand side, we can write 
\begin{eqnarray}
&& [U_{\alpha}(k-1,n) - L_{\alpha}(k-1,n)]\ + \Big{|} U_{\alpha}(k-1,n) - \left(\frac{k}{n+1}\right)^{\alpha}\Big{|}  \nonumber \\
&\leq& 2\left[\left(\frac{n+1+\alpha}{n+1}\right)^{1-\alpha} - 1\right] +  \left(\frac{k}{n+1}\right)^{\alpha}\left[1 - \left(\frac{k}{k+\alpha}\right)^{1-\alpha}\right] \ .  \label{eq:UnLnbis} 
\end{eqnarray}
At this stage, we observe that $\left(\frac{n+1+\alpha}{n+1}\right)^{1-\alpha} - 1 \leq \frac{\alpha(1-\alpha)}{n+1}$ holds for all $n \in \naturals$, while for the latter summand on the right-hand side of \eqref{eq:UnLnbis} we get
$$
\left(\frac{k}{n+1}\right)^{\alpha}\left[1 - \left(\frac{k}{k+\alpha}\right)^{1-\alpha}\right] \leq \alpha(1-\alpha)(1+\alpha)^{\alpha}\left(\frac{1}{n+1}\right)^{\alpha}\ .
$$
Moreover, for the third summand on the right-hand side of \eqref{eq:UnLn} we have $\left(\frac{k}{n}\right)^{\alpha} - \left(\frac{k}{n+1}\right)^{\alpha} \leq \frac{\alpha}{n}$, while for the last summand on the right-hand side of the same relation
we obtain $\left(\frac{k}{n}\right)^{\alpha} - \left(\frac{k-1}{n}\right)^{\alpha} \leq \left(\frac{1}{n}\right)^{\alpha}$. Putting these bounds together via \eqref{eq:UnLn}-\eqref{eq:UnLnbis}, we get \eqref{eq:beta} for $\alpha \in (0,1)$. When $\alpha \in (1,2)$, we start again from \eqref{eq:fellergamma}, which can be equivalently rewritten as
$$
\frac{\Gamma(n+1)}{\Gamma(n+1+\delta)} \frac{\Gamma(\lfloor nx\rfloor+1+\delta)}{\Gamma(\lfloor nx\rfloor+1)} \frac{\lfloor nx\rfloor+1+\delta}{n+1+\delta}
$$
with $\delta := \alpha -1$. Whence,
$$
\frac{\lfloor nx\rfloor+\alpha}{n+\alpha} L_{\delta}(\lfloor nx\rfloor,n) \leq \ff_n(x) \leq \frac{\lfloor nx\rfloor+\alpha}{n+\alpha} U_{\delta}(\lfloor nx\rfloor,n)
$$
for all $n \in \naturals$ and $x \in (0,1)$. Then, we can establish a bound similar to \eqref{eq:UnLn}, namely
\begin{align}
&\sup_{x \in [\frac{k-1}{n}, \frac{k}{n})} \Big{|} \ff_n\left(\frac{k-1}{n}\right) - x^{\alpha} \Big{|} \nonumber \\
&\quad\leq \frac{k+\delta}{n+\alpha}[U_{\delta}(k-1,n) - L_{\delta}(k-1,n)] + \frac{k+\delta}{n+\alpha}\Big{|} U_{\delta}(k-1,n) - \left(\frac{k}{n+1}\right)^{\delta} \Big{|} \nonumber \\
&\quad\quad+ \frac{k+\delta}{n+\alpha}\left[ \left(\frac{k}{n}\right)^{\delta} - \left(\frac{k}{n+1}\right)^{\delta} \right] + \left(\frac{k}{n}\right)^{\delta} \Big{|} \frac{k+\delta}{n+\alpha} - \frac{k}{n} \Big{|} + \left[ \left(\frac{k}{n}\right)^{\alpha} - \left(\frac{k-1}{n}\right)^{\alpha} \right]\ . \label{eq:UnLnter} 
\end{align}
We analyze the first two summands on the above right-hand side by resorting to \eqref{eq:UnLnbis}, to obtain 
\begin{align*}
&\frac{k+\delta}{n+\alpha}[U_{\delta}(k-1,n) - L_{\delta}(k-1,n)] + \frac{k+\delta}{n+\alpha}\Big{|} U_{\delta}(k-1,n) - \left(\frac{k}{n+1}\right)^{\delta} \Big{|} \nonumber \\
&\quad\leq 2\left[\left(\frac{n+1+\delta}{n+1}\right)^{1-\delta} - 1\right] +  2\left(\frac{k}{n+1}\right)^{\alpha}\left[1 - \left(\frac{k}{k+\delta}\right)^{1-\delta}\right] \nonumber\ .
\end{align*}
Now, the former summand on the above right-hand side has been already bounded by $\frac{2\alpha(1-\alpha)}{n+1}$, so that we can focus the attention on the latter. Arguing as above, we get
$$
2\left(\frac{k}{n+1}\right)^{\alpha}\left[1 - \left(\frac{k}{k+\delta}\right)^{1-\delta}\right] \leq 2\delta(1-\delta)(1+\delta)^{\delta} \frac{1}{n+1}\ .
$$
Lastly, the very same arguments used to handle the case $\alpha \in (0,1)$ lead to conclude that also the last three terms on the right-hand side of \eqref{eq:UnLnter} are bounded from above by a term of the type $C_{\alpha}/n$ for some constant $C_{\alpha}$ independent of $k$. The proof of \eqref{eq:beta} is therefore complete in the case that the prior is a beta distribution with parameters $(\alpha, 1)$.

The case of a beta distribution with parameters $(1,\alpha)$ is easily reformulated in terms a beta distribution with parameters $(\alpha, 1)$, in view of the following symmetry argument. First, we note that the d.f. $\ff(x)$ of a beta with parameters $(1,\alpha)$ coincides with $1- \ff^{\ast}(1-x)$, where $\ff^{\ast}$ is the d.f. of a beta with parameters $(\alpha,1)$. An analogous argument is true for the d.f. $\ff_n(x)$, in the sense that it coincides, for all $x \in [0,1]\setminus\{0, \frac{1}{n}, \frac{2}{n}, \dots, 1\}$, with the d.f. $\ff_n^{\ast}(x)$ of the r.v. $\frac{1}{n} \sum_{i=1}^n \overline{X}_i$, where $\overline{X}_i := 1-X_i$ for all $i \in \naturals$. The conclusion is reached by observing that the de Finetti measure of the (exchangeable) sequence $\{\overline{X}_i\}_{i \geq 1}$ is exactly the beta with parameters $(\alpha, 1)$, whenever the de Finetti measure of the sequence $\{X_i\}_{i \geq 1}$ is the beta with parameters $(1,\alpha)$, and that
\begin{align*}
\ud_K(\mu_n; \mu) &= \sup_{x \in [0,1]\setminus\{0, \frac{1}{n}, \frac{2}{n}, \dots, 1\}} |\ff_n(x) - \ff(x)|  = 
\sup_{x \in [0,1]\setminus\{0, \frac{1}{n}, \frac{2}{n}, \dots, 1\}} |\ff_n^{\ast}(x) - \ff^{\ast}(x)| \\
&= \sup_{x \in [0,1]} |\ff_n^{\ast}(x) - \ff^{\ast}(x)| \ . \nonumber 
\end{align*}

\subsection{Proof of Proposition \ref{prop:Holder}}

By resorting once again to formula \eqref{eq:feller} in the proof of Proposition \ref{prop:beta}, we get
\begin{eqnarray}
\ud_K(\mu_n; \mu) &\leq& \sup_{x \in [0,1)} \int^1_0 \beta(\theta; \lfloor nx\rfloor + 1, n - \lfloor nx\rfloor) |\ff(x) - \ff(\theta)| \ud \theta \nonumber \\
&\leq& H_{\gamma}(\ff) \sup_{x \in [0,1)} \int^1_0 \beta(\theta; \lfloor nx\rfloor + 1, n - \lfloor nx\rfloor) |x - \theta|^{\gamma} \ud\theta \nonumber
\end{eqnarray}
where $H_{\gamma}(\ff)$ denotes the H\"older constant of $\ff$. Now, we exploit that $|x - \theta|^{\gamma} \leq |x - \eta_x|^{\gamma} + |\eta_x - \theta|^{\gamma}$, where 
$\eta_x :=
\frac{\lfloor nx\rfloor + 1}{n+1} = \int^1_0 \theta \beta(\theta; \lfloor nx\rfloor + 1, n - \lfloor nx\rfloor ) \ud\theta$, to get 
$$
\ud_K(\mu_n; \mu) \leq H_{\gamma}(\ff) \Big[ \sup_{x \in [0,1)} |x - \eta_x|^{\gamma} + \sup_{x \in [0,1)} \int^1_0 |\eta_x - \theta|^{\gamma} \beta(\theta; \lfloor nx\rfloor + 1, n - \lfloor nx\rfloor) \ud\theta \Big]\ .
$$
Since $\sup_{x \in [0,1)} |x - \eta_x| \leq \frac{2}{n+1}$ follows from direct computation, we can focus on the second summand on the above right-hand side, which can be bounded by means of the Jensen inequality as follows:
$$
\int^1_0 |\eta_x - \theta|^{\gamma} \beta(\theta; \lfloor nx\rfloor + 1, n - \lfloor nx\rfloor) \ud\theta \leq \Big(\int^1_0 |\eta_x - \theta|^2 \beta(\theta; \lfloor nx\rfloor  + 1, n - \lfloor nx\rfloor) \ud\theta\Big)^{\gamma/2}\ .
$$
The proof is completed by observing that the integral $\int^1_0 |\eta_x - \theta|^2 \beta(\theta; \lfloor nx\rfloor  + 1, n - \lfloor nx\rfloor) \ud\theta$ represents the variance of the beta distribution with parameters
$(\lfloor nx\rfloor  + 1, n - \lfloor nx\rfloor)$, which, being equal to $\frac{(\lfloor nx\rfloor+1)(n-\lfloor nx\rfloor )}{(n+1)^2(n+2)}$, is less than $\frac{1}{n+2}$ for any $x \in [0,1]$.

\section*{Acknowledgements}

The authors thank two anonymous Referees for all their comments, corrections, and suggestions which remarkably improved the original version of the paper. Emanuele Dolera and Stefano Favaro received funding from the European Research Council (ERC) under the European Union's Horizon 2020 research and innovation programme under grant agreement No 817257. Emanuele Dolera and Stefano Favaro gratefully acknowledge the financial support from the Italian Ministry of Education, University and Research (MIUR), ``Dipartimenti di Eccellenza" grant 2018-2022.


\begin{flushleft}

EMANUELE DOLERA\\  
Universit\`a di Pavia\\ 
Dipartimento di Matematica\\
Via Adolfo Ferrata 5, 27100 Pavia, Italy\\
emanuele.dolera@unipv.it

\bigskip

%
STEFANO FAVARO\\
Universit\`a  degli studi di Torino\\ 
Dipartimento di Economia e Statistica\\ 
Lungo Dora Siena 100A, 10134 Torino, Italy\\ 
stefano.favaro@unito.it

\end{flushleft}

\end{document}